\documentclass[10pt]{article}
\usepackage{flafter,mathrsfs,amsmath,amssymb,latexsym,psfrag,graphicx,color}
\usepackage[margin=2.5cm]{geometry}
\usepackage[numbers,sort&compress]{natbib}
\allowdisplaybreaks

\newtheorem{theorem}{Theorem}[section]
\newtheorem{defn}[theorem]{Definition}
\newtheorem{lemma}[theorem]{Lemma}
\newtheorem{remark}[theorem]{Remark}

\newtheorem{claim}[theorem]{Claim}
\makeatletter

\renewcommand{\theequation}{\thesection.\arabic{equation}}

\@addtoreset{figure}{section}
\makeatother

\begin{document}
\setlength\arraycolsep{2pt}
\date{\today}

\title{Solvability of interior transmission problem for the diffusion equation by constructing its Green function}
\author{Gen Nakamura$^1$, Haibing Wang$^{2}$
\\$^1$Department of Mathematics, Hokkaido University, Sapporo 060-0810, Japan
\\E-mail: nakamuragenn@gmail.com
\\$^2$School of Mathematics, Southeast University, Nanjing 210096, P.R. China \,\quad
\\E-mail: hbwang@seu.edu.cn
}

\maketitle

\begin{abstract}

Consider the interior transmission problem arising in inverse boundary value problems for the diffusion equation with discontinuous diffusion coefficients. We prove the unique solvability of the interior transmission problem by constructing its Green function. First, we construct a local parametrix for the interior transmission problem near the boundary in the Laplace domain, by using the theory of pseudo-differential operators with a large parameter. Second, by carefully analyzing the analyticity of the local parametrix in the Laplace domain and estimating it there, a local parametrix for the original parabolic interior transmission problem is obtained via the inverse Laplace transform. Finally, using a partition of unity, we patch all the local parametrices and the fundamental solution of the diffusion equation to generate a global parametrix for the parabolic interior transmission problem, and then compensate it to get the Green function by the Levi method. The uniqueness of the Green function is justified by using the duality argument, and then the unique solvability of the interior transmission problem is concluded. We would like to emphasize that the Green function for the parabolic interior transmission problem is constructed for the first time in this paper. It can be applied for active thermography and diffuse optical tomography modeled by diffusion equations to identify an unknown inclusion and its physical property.

\bigskip



\end{abstract}

\section{Introduction}
\setcounter{equation}{0}

Interior transmission problem plays an important role in inverse scattering theory for inhomogeneous media. It is a non-classical boundary value problem for a pair of partial differential equations in a bounded domain coupled on the boundary. There are many works on interior transmission problems for elliptic equations; see, for example, \cite{C-C-D2010, C-C-H2010, C-C-H2016, C-G-H2010, C-H2011, C-P-S2007, Fai2014, G-V-Z2016, H-K-O-P2010, H-K-O-P2011, Kir1986, Kir2009, L-V2012, L-V2015, L-R2015, M-S2012, P-S2008, R-S1991, Sun2011, Syl2012, Syl2013}, studying the solvability of interior transmission problems for different kinds of inhomogeneous media, the existence and efficient computations of the transmission eigenvalues, and their applications to inverse scattering problems. Recently, it was found in \cite{J-L2017, L-L-L-W2017} that interior transmission problems are closely related to the invisibility cloak in acoustic and electromagnetic wave scattering.

As we explained in our paper \cite{N-W2013}, the interior transmission problem for the diffusion equation naturally arises in inverse boundary value problems for the diffusion equation with discontinuous coefficients when we consider a reconstruction method called the linear sampling method. More explicitly, the problem arises from studying the solvability of the so-called Neumann-to-Dirichlet map (ND map) equation, which plays a central role in the linear sampling method for reconstructing unknown inclusions inside a diffusive medium from boundary measurements. Let $D$ be an inclusion compactly embedded in the diffusive medium. Assume that the diffusion coefficients of $D$ and the background are $k$ and $1$, respectively. We assume for simplicity that $k$ is a constant with $k>1$. Suppose that $D\subset\mathbb{R}^n\,(n=2, 3)$ is a bounded domain with $C^\infty$ smooth boundary $\partial D$. Then the interior transmission problem (ITP) for the diffusion equation is described by the following initial boundary value problem:
\begin{equation}\label{eq:ITP1}
\begin{cases}
(\partial_t-\Delta)v=0 & \textrm{ in }D_T:=D\times(0,T),\\
(\partial_t-k\Delta)u=0 &  \textrm{ in }D_T:=D\times(0,T),\\
v-u=f & \textrm{ on }(\partial D)_T:=\partial D\times (0,T),\\
\partial_{\nu}v - k\partial_{\nu}u=g & \textrm{ on }(\partial D)_T:=\partial D\times (0,T),\\
v=u=0 & \textrm{ at } t=0,
\end{cases}
\end{equation}
where $\nu$ is the unit outer normal vector to $\partial D$. In general, the interior transmission problem may have inhomogeneous terms in the equations.

One may think that the interior transmission problem is a special problem attached to the linear sampling method for the aforementioned inverse boundary value problem. We would like to show by giving a general example that this is not true. Let us consider $D$ as a heat conductor with thermal conductivity $c$ located in an open air with thermal conductivity $c_0$. Put some heat source $p(x,\,t)$ over the time interval $(0,\,T)$ and let it radiate. Then the temperature $u=u(x,\,t)$ generated by this heat source satisfies the following initial value problem:
\begin{equation}
\begin{cases}
\partial_t u - \nabla\cdot(\gamma\nabla u)=p & \textrm{ in } {\mathbb R}^3\times(0,T),\\
u=0 & \textrm{ at } t=0,
\end{cases}
\end{equation}
where $\gamma=c_0+(c-c_0)\chi_D$ with $\chi_D$ being the characteristic function of $D$. Suppose we can measure
\begin{equation}\label{Meas}
u\big|_{(\partial D)_T},\quad \partial_{\gamma}u:=\nu\cdot(c\nabla u)\big|_{(\partial D)_T}=c\partial_\nu u\big|_{(\partial D)_T},
\end{equation}
and consider the problem of identifying the unknown conductivity $c$ from this measurement. That is, for two unknowns $c_1$ and $c_2$, show that they are equal if we have $u_1=u_2$ and $c_1\partial_\nu u_1=c_2\partial_\nu u_2$ on $(\partial D)_T$, where $u_j\,(j=1,\,2)$ satisfies
\begin{equation}\label{itpRadiation}
\begin{cases}
(\partial_t-c_j\Delta)u_j=p & \textrm{ in } D_T,\\
u_j=0 & \textrm{ at } t=0.
\end{cases}
\end{equation}
Combining \eqref{itpRadiation} with the measurement \eqref{Meas} for $u=u_j\,(j=1,\,2)$, we have an interior transmission problem for $(u_1,\,u_2)$ in $D_T$.

In \eqref{eq:ITP1}, we assume that $f=g=0$ at $t=0$ and $f,\,g$ satisfy certain regularity assumption; say $f\in H^1((0,\,T);\,H^{3/2}(\partial D))$ and $g\in H^1((0,\,T);\,H^{1/2}(\partial D))$. Then we can remove $f$ and $g$ so that we have homogeneous boundary conditions and inhomogeneous terms in the equations. So we are led to the following parabolic interior transmission problem:
\begin{equation}\label{eq:ITP-2}
\begin{cases}
(\partial_t-\Delta)v=N_1 & \textrm{ in }D_T,\\
(\partial_t-k\Delta)u=N_2 & \textrm{ in }D_T,\\
v-u=0 & \textrm{ on }(\partial D)_T,\\
\partial_{\nu}v -k\partial_{\nu}u=0 & \textrm{ on }(\partial D)_T, \\
v=u=0 & \textrm{ at } t=0.
\end{cases}
\end{equation}

In this paper, we will show the unique solvability of \eqref{eq:ITP-2} by constructing its Green function; see Section \ref{patch} for the definition of the meaning of the unique solvability of \eqref{eq:ITP-2}. In order to define the Green function $\mathbb{G}$ for \eqref{eq:ITP-2}, let
\begin{equation}\label{F ell}
F_\ell:=\big((2-\ell)\delta(x-y,\,t-s),\, (\ell-1)\delta(x-y,\,t-s)\big) \qquad \textrm{with}\; \,y\in D,\,s\in(0,\,T)
\end{equation}
for $\ell=1,\,2$. Define the distribution $K_\ell(x,\,t;\,y,\,s)=\big(G_\ell(x,\,t;\,y,\,s),\, H_\ell(x,\,t;\,y,\,s)\big)\in\mathscr D^\prime(D_T\times D_T)$ as the solution of the initial boundary value problem \eqref{eq:ITP-2} with $(N_1,\,N_2)=F_\ell$ such that $K_\ell(x,\,t;\,y,\,s)=(0,\,0)$ for $0<t<s$ when we fix $(y,\,s)\in D_T$. Then the Green function is expressed by the matrix
\begin{equation}\label{Green function G}
\mathbb G := \left((K_1)^{\mathbb T},\,(K_2)^{\mathbb T}\right)=\left(
\begin{array}{lr}
G_1 & G_2 \\
H_1 & H_2
\end{array} \right),
\end{equation}
where \lq\lq $\mathbb T$\rq\rq{} denotes the transpose. Based on this observation, we perform the construction of the Green function in the following way. First, we construct a local parametrix for \eqref{eq:ITP-2} near the boundary, by studying \eqref{eq:ITP-2} with $(N_1,\,N_2)=F_\ell$. To this end, we consider the interior transmission problem in the Laplace domain
\begin{equation}\label{eq:L_ITP-2}
\begin{cases}
(\tau-\Delta)G_\ell=(2-\ell)e^{-\tau s}\delta(x-y) & \mbox{ in }D,\\
(\tau-k\Delta)H_\ell=(\ell-1)e^{-\tau s}\delta(x-y) & \mbox{ in }D,\\
G_\ell-H_\ell=0 & \mbox{ on }\partial D,\\
\partial_{\nu}G_\ell-k\partial_{\nu}H_\ell=0 & \mbox{ on }\partial D.
\end{cases}
\end{equation}
Here $\tau\in{\mathbb C}$ denotes the Laplace variable in a sector of $\mathbb C$ where $|\tau|$ can become large (see Definition \ref{Def-sym} and Section \ref{estimate}). By using the theory of pseudo-differential operators with a large parameter \cite{Shubin}, we construct a local parametrix for \eqref{eq:L_ITP-2} in a neighborhood of $x_0\in \partial D$, which is analytic with respect to $\tau$ in some restricted domain and has a good estimate. This leads us to a local parametrix for \eqref{eq:ITP-2} near the boundary via the inverse Laplace transform. Second, using a partition of unity, we patch all the local parametrices and the fundamental solution of the diffusion equation to generate a global parametrix for \eqref{eq:ITP-2}, and then compensate it to obtain the desired Green function by the Levi method. Finally, we show the uniqueness of the Green function by the duality argument.

The novelty and new contributions of the present work are as follows. First, by showing the solvability of the interior transmission problem \eqref{eq:ITP-2} as a by-product of the construction of its Green function, we could clearly clarify the solvability of the ND map equation, which strengthens our theoretical analysis on the sampling method proposed in \cite{N-W2013}. Second, our argument of constructing the Green function for \eqref{eq:ITP-2} is new and it is efficient to see its principal part very clearly. The argument is a adaptation of Seeley's argument \cite{Seeley} for elliptic boundary value problems to interior transmission problems for diffusion equations. In \cite{K-N-S} one of the authors of this paper showed how to adapt Seeley's argument to construct the Green function for the elliptic interior transmission problem assuming its unique solvability. Concerning the construction of a local parametrix for the Green function, the argument in \cite{K-N-S} is much more simple than the direct application of Seeley's argument. Since in our case we are studying the unique solvability of \eqref{eq:ITP-2} via constructing its Green function, we have to compensate a parametrix to obtain the Green function without using the unique solvability. We could achieve it by using the Levi method. Here we remark that the solvability of \eqref{eq:ITP-2} is the consequence of the existence of its Green function and it is not available beforehand. In addition, our argument gives the leading part of the Green function very clearly, which can be utilized to analyze the asymptotic behavior of the indicator function of the linear sampling method which is a non-iterative reconstruction method for inverse boundary value problems for parabolic equations, since the indicator function is closely related to the corresponding Green function \cite{I-K-N2010, K-N-S, N-W2015, N-W2017, W-L2018}. Let us show the background about this in more details.

Let $\Omega$ be a thermal conductor and $D$ an inclusion inside the conductor. The thermal conductivities of $\Omega$ and $D$ are $1$ and $k$, respectively. For any given heat flux $g$, the temperature $u(x,\,t)$ satisfies the following initial-boundary value problem:
\begin{equation}\label{eq:mp}
\begin{cases}
(\partial_t - \nabla \cdot k \nabla)u=0 & \textrm{ in } D_T, \\
(\partial_t - \Delta)u=0 & \textrm{ in } (\Omega\setminus\overline D)_T, \\
u|_- - u|_+=0 & \textrm{ on } (\partial D)_T, \\
k\partial_\nu u|_--\partial_\nu u|_+=0 & \textrm{ on } (\partial D)_T, \\
\partial_\nu u=g & \textrm{ on } (\partial \Omega)_T,\\
u=0 & \textrm{ at } t=0,
\end{cases}
\end{equation}
where $\nu$ on $\partial D$ (or $\partial \Omega$) is the unit normal vector directed into the exterior of $D$ (or $\Omega$). Here the subscripts \lq\lq$+$\rq\rq{} and \lq\lq$-$\rq\rq{} indicate the trace taken from the exterior and interior of $D$, respectively. The inverse problem is to identify the unknown inclusion $D$ from the Neumann-to-Dirichlet map $\Lambda_D$. As studied in \cite{N-W2013, W-L2018}, we can characterize the boundary of $D$ by solving the Neumann-to-Dirichlet map gap equation
\begin{equation}\label{gap0}
(\Lambda_D-\Lambda_\emptyset)g=G^\Omega_{(y,\,s)}(x,\,t),
\end{equation}
where $\Lambda_\emptyset$ is the Neumann-to-Dirichlet map when $D=\emptyset$, and $G^\Omega_{(y,\,s)}(x,\,t):=G^\Omega(x,\,t;\,y,\,s)$ is the Green function for the heat operator $\partial_t - \Delta$ in $\Omega_T$ with homogeneous Neumann boundary condition on $(\partial \Omega)_T$. Then, for $y\in D$ and $s\in (0,\,T)$, the equation \eqref{gap0} has a solution if and only if the interior transmission problem
\begin{equation}\label{ITP}
\begin{cases}
(\partial_t-\Delta)v=0 & \mathrm{in}\;D_T,\\
(\partial_t - \nabla\cdot k\nabla)w=0 & \mathrm{in}\;D_T,\\
w-v=G^\Omega_{(y,\,s)}(x,\,t) & \mathrm{on}\;(\partial D)_T,\\
k\partial_\nu w-\partial_\nu v=\partial_\nu G^\Omega_{(y,\,s)}(x,\,t) &  \mathrm{on}\;(\partial D)_T,\\
w=0 & \mathrm{at}\;t=0,\\
v=0 & \mathrm{at}\;t=0
\end{cases}
\end{equation}
is solvable with the solution $w$ and $v$ satisfying the equations $(\partial_t-\nabla\cdot\gamma\nabla)w=0$ and $(\partial_t-\Delta)v=0$ in $\Omega_T$, respectively, where $\gamma = 1+(k-1)\chi_D$. Moreover, if $g$ is the solution to \eqref{gap0}, we let $v$ satisfy $(\partial_t - \Delta )v=0$ in $\Omega_T$ with zero initial condition and the boundary data $\partial_\nu v|_{(\partial \Omega)_T}=g$. Then we have
\begin{equation}\label{psol5}
v = G^D_{(y,\,s)}(x,\,t) - G^\Omega_{(y,\,s)}(x,\,t) \quad \mathrm{in}\; D_T,
\end{equation}
where $G^D_{(y,\,s)}(x,\,t)$ meets
\begin{equation}\label{psol6}
\begin{cases}
(\partial_t - \nabla\cdot k\nabla) H^D_{(y,\,s)}=0 & \mathrm{in}\;D_T,\\
(\partial_t-\Delta) G^D_{(y,\,s)}=\delta(x-y)\,\delta(t-s) & \mathrm{in}\;D_T,\\
H^D_{(y,\,s)} - G^D_{(y,\,s)} = 0 & \mathrm{on}\;(\partial D)_T,\\
k\partial_\nu H^D_{(y,\,s)} - \partial_\nu G^D_{(y,\,s)} = 0 &  \mathrm{on}\;(\partial D)_T,\\
H^D_{(y,\,s)}=0 & \mathrm{at}\;t=0,\\
G^D_{(y,\,s)}=0 & \mathrm{at}\;t=0.
\end{cases}
\end{equation}
We clearly see that the solution $(G^D_{(y,\,s)},\,H^D_{(y,\,s)})$ to \eqref{psol6} is the second column of the Green matrix \eqref{Green function G}. For more details, we can see the reference \cite{W-L2018}.

The paper is organized as follows. In Section \ref{flatP}, we construct a local parametrix for the interior transmission problem in the Laplace domain by studying \eqref{eq:L_ITP-2}. Some lengthy details appear in the appendix. Then, in Section \ref{estimate}, taking the inverse Laplace transform of this parametrix, we obtain a local parametrix for \eqref{eq:ITP-2} with some estimates. In Section \ref{patch}, using a partition of unity, we patch the local parametrices and the fundamental solution of the diffusion equation so that we have a global parametrix for \eqref{eq:ITP-2}. This parametrix can be compensated to generate the Green function by the Levi method. In Section \ref{unique}, the uniqueness of the Green function is justified, and the unique solvability of \eqref{eq:ITP-2} is summarized. Finally, in Section \ref{conclusion}, we give some concluding remarks.

\section{Construction of a local parametrix in the Laplace domain}\label{flatP}
\setcounter{equation}{0}

In this section, we construct a local parametrix for the interior transmission problem in the Laplace domain by studying \eqref{eq:L_ITP-2}. Our argument is based on the theory of pseudo-differential operators with a large parameter. We only consider the case of $n=3$. The construction of the local parametrix is performed as follows. First, we locally flatten the boundary $\partial D$ by a coordinate transformation, and transform \eqref{eq:L_ITP-2} locally into an elliptic system \eqref{eq:L_ITP-3} defined in a half space. Then, for solving \eqref{eq:L_ITP-3}, we equivalently solve the transmission problem \eqref{eqns}--\eqref{boundc}. The solution $(G^\pm_\ell,\,H^\pm_\ell)$ to this transmission problem is obtained by regarding the functions $G^\pm_\ell$ and $H^\pm_\ell$ as Schwartz kernels of the corresponding pseudo-differential operators. Thus, we are led to determine the amplitudes of the pseudo-differential operators, which is accomplished in Theorem \ref{param}.

Let us first locally flatten the boundary $\partial D$ near a point $x_0\in \partial D$ by a coordinate transformation $\Phi_{x_0}:\, U(x_0)\to \mathbb R^3$ with $\Phi_{x_0}(x_0)=0$, where $U(x_0)$ is an open neighbourhood of $x_0$ in $\mathbb R^3$. Under this coordinate transformation, we can locally express $\partial D$ and $D$ by $\partial D = \{\xi_3=0 \}$ and $D = \{\xi_3<0 \}$, respectively, in terms of the local coordinates $\xi=(\xi_1,\,\xi_2,\,\xi_3)$. Denote by $J:=\nabla_x \xi$ the Jacobian of the coordinate transformation. Define $M=\left( m_{jl} \right)_{3\times 3}:=J J^T$ and $\mathcal J := \mathrm{det}\,(\nabla_\xi x)$. Without loss of generality, we assume $\mathcal J>0$ by fixing the orientation of $\partial D$ considered as a manifold. Let $\xi=\Phi_{x_0}(x),\,\eta=\Phi_{x_0}(y)$, $G^\sharp_\ell (\xi)=G_\ell(\Phi^{-1}_{x_0}(\xi)),\,H^\sharp_\ell (\xi)=H_\ell(\Phi^{-1}_{x_0}(\xi))$. Then from \eqref{eq:L_ITP-2} we locally have
\begin{equation}\label{eq:L_ITP-3}
\begin{cases}
\mathcal PG^\sharp_\ell := \tau G^\sharp_\ell - \mathcal J^{-1}(\xi)\nabla_\xi \cdot (\mathcal J(\xi)M\nabla_\xi G^\sharp_\ell)= (2-\ell) e^{-\tau s}\mathcal J^{-1}(\eta)\delta(\xi-\eta) & \mathrm{in}\;\mathbb R^3_{-},\\
\mathcal QH^\sharp_\ell :=\tau H^\sharp_\ell - k\mathcal J^{-1}(\xi)\nabla_\xi \cdot (\mathcal J(\xi)M\nabla_\xi H^\sharp_\ell)= (\ell -1) e^{-\tau s}\mathcal J^{-1}(\eta)\delta(\xi-\eta) & \mathrm{in}\;\mathbb R^3_{-},\\
G^\sharp_\ell - H^\sharp_\ell=0 & \mathrm{on}\;\partial \mathbb R^3_{-},\\
e_3\cdot M\nabla_\xi G^\sharp_\ell - ke_3\cdot M\nabla_\xi H^\sharp_\ell = 0 & \mathrm{on}\;\partial \mathbb R^3_{-},
\end{cases}
\end{equation}
where $\mathbb R^3_-:=\{(\xi_1,\,\xi_2,\,\xi_3)\in \mathbb R^3:\, \xi_3<0 \}$ and $\partial \mathbb R^3_-:=\{(\xi_1,\,\xi_2,\,\xi_3)\in \mathbb R^3:\, \xi_3=0 \}$.
$\mathcal P$ and $\mathcal Q$ are strongly elliptic second order operators on $\Phi_{x_0}(U(x_0))\cap \overline{\mathbb R_-^3}$ with smooth coefficients. We extend them to the whole $\overline{\mathbb R_-^3}$ without destroying their strong ellipticity and the smoothness of their coefficients. Hence, we will look for $(G^\sharp_\ell,\,H^\sharp_\ell)$ which satisfies \eqref{eq:L_ITP-3} in the whole $\overline{\mathbb R_-^3}$.

In the sequel, for convenience, we will still use the notations $G_\ell,\,H_\ell,\,x,\,y$ in the local coordinates system, instead of $G^\sharp_\ell,\,H^\sharp_\ell,\,\xi,\,\eta$. Let $y=(y_1,\,y_2,\,y_3)\in \mathbb R^3_-$ with $y_3$ near to $0$. To clarify the dependency of $G_\ell$ and $H_\ell$ on $y$, we denote them by $G_\ell(x,\,y)$ and $H_\ell(x,\,y)$, respectively, where we have suppressed $\tau$. We further represent them in the forms:
\begin{equation}\label{G_ell}
G_\ell(x,\,y)=\begin{cases}
G_\ell^{+}(x,\,y), & \textrm{ if }x_{3}-y_{3}>0,\\
G_\ell^{-}(x,\,y), & \textrm{ if }x_{3}-y_{3}<0,
\end{cases}
\end{equation}
and
\begin{equation}\label{H_ell}
H_\ell(x,\,y)=\begin{cases}
H_\ell^{+}(x,\,y), & \textrm{ if }x_{3}-y_{3}>0,\\
H_\ell^{-}(x,\,y), & \textrm{ if }x_{3}-y_{3}<0.
\end{cases}
\end{equation}
Then $G_\ell^{\pm}$ and $H_\ell^{\pm}$ satisfy the equations
\begin{equation}\label{eqns}
\mathcal PG_\ell^{\pm}=\mathcal QH_\ell^{\pm}=0 \quad \textrm{ for }\pm(x_{3}-y_{3})>0,\,x_{3}<0,
\end{equation}
the transmission conditions on $x_3=y_3$
\begin{equation}\label{tranc}
\begin{cases}
G_\ell^{+}=G_\ell^{-},\\
H_\ell^{+}=H_\ell^{-},\\
e_3\cdot M(x)\nabla(G_\ell^{+}-G_\ell^{-})=-(2-\ell)e^{-\tau s} \mathcal J^{-1}(y)\, \delta(x'-y'),\\
ke_3\cdot M(x)\nabla(H_\ell^{+}-H_\ell^{-})=-(\ell-1)e^{-\tau s} \mathcal J^{-1}(y)\, \delta(x'-y')
\end{cases}
\end{equation}
with $x'=(x_1,\,x_{2}),\,y'=(y_1,\,y_2)$, and the boundary conditions on $x_3=0$
\begin{equation}\label{boundc}
\begin{cases}
G_\ell^{+}=H_\ell^{+},\\
e_3\cdot M(x)\nabla G_\ell^{+}=ke_3\cdot M(x)\nabla H_\ell^{+}.
\end{cases}
\end{equation}
We note that solving \eqref{eq:L_ITP-3} is equivalent to solving \eqref{G_ell}--\eqref{boundc}. Also, the transmission conditions on $x_3=y_3$ are coming from the two equations in \eqref{eq:L_ITP-3}. The functions $G_\ell^\pm$ and $H_\ell^\pm$ will be obtained as Schwartz kernels of pseudo-differential operators with symbols in $S(\infty)$ which is given as follows.

\begin{defn} \label{Def-sym}
For $m\in \mathbb R$, $\mathbf{a}(x^\prime,\,\xi^\prime,\,\tau)$ is in $S(m)$ if the followings hold:
\begin{itemize}
\item[$(1)$] $\mathbf{a}(x^\prime,\,\xi',\,\tau)\in C^{\infty}(\mathbb R^2_{x'} \times \mathbb R^2_{\xi'} \times \Sigma(\xi^\prime))$;

\item[$(2)$] For arbitrary multi-indices $\alpha,\, \beta \in \mathbb Z^2_{+}$ with $\mathbb Z_+:=\mathbb N\cup \{0\}$, there exists a constant $C_{\alpha,\,\beta}>0$ such that
\begin{equation*}
|D_{x'}^{\alpha}D_{\xi'}^{\beta}\mathbf a(x^\prime,\,\xi',\,\tau)|\leq C_{\alpha,\,\beta} \langle \xi^\prime,\,\tau \rangle^{m-|\beta|}, \quad
(x',\,\xi',\,\tau)\in \mathbb R^2_{x'}\times \mathbb R^2_{\xi'} \times \{\Sigma(\xi^\prime)\cap\{|\tau|>1\}\},
\end{equation*}
\end{itemize}
where $\langle \xi^\prime,\,\tau \rangle^2 = 1 + |\xi^\prime|^2 + |\tau|$ and $\Sigma(\xi^\prime):=\{ \sigma + re^{i\theta}: \, r>0,\,\theta_1 < \theta < \theta_2 \}$ with some constants $\sigma, ,\theta_1,\,\theta_2$ depending on $\xi^\prime$.

We call $\mathbf a$ a symbol of order $m$. Further, we define $S(\infty)$ by $S(\infty)=\bigcup _{m\in \mathbb R}S(m)$.

If $\mathbf a(x^\prime,\,\xi^\prime,\,\tau)=\mathbf a(x^\prime,\,\xi^\prime,\,\tau;\,x_3,\,y,\,s)$ depends on the parameters $x_3\in\mathbb R_-:=\{x_3:\,x_3<0\}$, $y\in {\mathbb R}_-^3:=\{y=(y^\prime,\,y_3)\in{\mathbb R}^3:\,y_3<0\}$ and $s\in (0,\,T)$, the above estimate has to be held uniformly with respect to them. For simplicity, we write $\mathbf a(x,\,\xi^\prime,\,\tau)=\mathbf a(x^\prime,\,\xi^\prime,\,\tau;\,x_3)$ and $\mathbf a(x,\,\xi^\prime,\,\tau;\,y,\,s)=\mathbf a(x^\prime,\,\xi^\prime,\,\tau;\,x_3,\,y,\,\,s)$.

For $\mathbf{a}(x^\prime,\,\xi^\prime,\,\tau)\in S(\infty)$, define the pseudo-differential operator $\mathbf{a}(x^\prime,\,D_{\xi^\prime},\,\tau)$ by
\begin{equation}\label{pdo}
\big(\mathrm{Op}(\mathbf a)\varphi\big)(x^\prime,\,\tau)=\big(\mathbf{a}(x^\prime,\,D_{\xi^\prime},\,\tau)\varphi\big)(x^\prime,\,\tau)=(2\pi)^{-2} \int_{\mathbb R^2} e^{ix^\prime \cdot \xi^\prime} \mathbf a(x^\prime,\,\xi^\prime,\,\tau)\, (\mathcal F  \varphi)(\xi^\prime)\,d\xi^\prime,
\end{equation}
where $$(\mathcal F  \varphi)(\xi^\prime)=\int_{{\mathbb R}^2}e^{-ix^\prime\cdot\xi^\prime}\varphi(x^\prime)\,dx^\prime$$
is the Fourier transform of $\varphi(x^\prime)\in C_0^\infty({\mathbb R}^2)$. Associated with $S(m)$ and $S(\infty)$, we define $S[m]:=\{\mathrm{Op}(\mathbf a):\,\mathbf a\in S(m)\}$ and $S[\infty]:=\bigcup _{m\in\mathbb R}S[m]$, respectively.
\end{defn}

Consider the operator
\begin{equation*}
\mathcal{P} = -\nabla \cdot (M(x)\nabla) - W \cdot \nabla + \tau = D\cdot(M(x)D) - i W\cdot D + \tau,
\end{equation*}
where $W=(w_1,\,w_2,\,w_3)^T:=\mathcal J^{-1} M(x) \nabla \mathcal J$ and $D=-i\nabla=(D_1,\,D_2,\,D_3)=:(D_{x^\prime},\,D_3)$. Rewrite $\mathcal P$ into the form
\begin{equation*}
\mathcal P = \sum_{m=0}^2 \mathcal P_m(x,\,D_{x^\prime},\,\tau) D_3^{2-m},
\end{equation*}
where each $\mathcal P_m$ is a partial differential operator of order $m$ with respect to $x^\prime$ depending on $x_3$ and $\tau$, and it can be considered as a pseudo-differential operator in $S[m]$. Hence, $\mathcal P$ can be viewed as a second order ordinary differential operator with respect to $x_3$ with coefficients in $S[\infty]$ and it is denoted by $\mathcal P=\mathcal P(x_3,\,D_3;\,x^\prime,\,D_{x^\prime},\,\tau)$. Decompose $\mathcal P$ into
\begin{equation*}
\mathcal P=p_2+p_1
\end{equation*}
with
\begin{align*}
& p_2(x_3,\,D_3;\,x^\prime,\,D_{x^\prime},\,\tau) = m_{33}(x)D_3^2 + 2\sum_{j=1}^2 m_{3j}(x)D_j D_3 + \sum_{j,l=1}^2 m_{jl}(x)D_j D_l + \tau,\\
& p_1(x_3,\,D_3;\,x^\prime,\,D_{x^\prime},\,\tau) = -\big( i\sum_{j=1}^3 \partial_{x_j}m_{j3}(x) + iw_3(x) \big) D_3 - \big( i
\sum_{j=1}^3\sum_{l=1}^2\partial_{x_j} m_{jl}(x)D_l  + i \sum_{j=1}^2 w_j(x)D_j \big),
\end{align*}
and expand each coefficient of $p_1,\,p_2$; say $p(x,\,D_{x^\prime},\,\tau)\in S[\infty]$, into its Taylor series around $x_3=y_3$. That is, expand $p(x,\,\xi',\,\tau)$ into
\begin{equation*}
p(x,\,\xi',\,\tau)=\sum_{j=0}^{\infty}(j!)^{-1}(x_3-y_3)^j (\partial_{x_3}^j
p)(x',\,y_3,\,\xi',\,\tau).
\end{equation*}

For our further arguments, we introduce the concept of order for the symbols in $S(\infty)$.
\begin{defn}\label{Def-order}
The multiplications by $\xi_1,\,\xi_2,\,\tau^{1/2}$ and $D_3$ are regarded as operators of order $1$, the multiplication by
$x_3-y_3$ is regarded as an operator of order $-1$, and $D_{x^\prime}$ is regarded as an operator of order $0$. The actual meaning of the order is as follows. If $\mathbf a(x,\,\xi^\prime,\,\tau)\in S(m)$ such that
$$D_{x^\prime}^\alpha D_{\xi^\prime}^\beta \mathbf a(x,\,\xi^\prime,\,\tau) \langle \xi^\prime,\,\tau \rangle^{-(m-|\beta|)} e^{\delta \langle \xi^\prime,\,\tau \rangle |x_3-y_3|} \qquad \textrm{with }\, x^\prime,\,\xi^\prime \in \mathbb R^2,\,x_3\leq 0,\,\tau\in \{\Sigma(\xi^\prime)\cap\{|\tau|>1\}\}$$
is bounded for small $\delta>0$ and each $\alpha,\,\beta\in \mathbb Z^2_+$, then we have for $j=1,\,2$ that
\begin{eqnarray*}
&&|D_{x^\prime}^\alpha D_{\xi^\prime}^\beta (\xi_j \mathbf a(x,\,\xi^\prime,\,\tau))|,\, |\tau^{1/2}D_{x^\prime}^\alpha D_{\xi^\prime}^\beta \mathbf a(x,\,\xi^\prime,\,\tau)|,\,|D_3 D_{x^\prime}^\alpha D_{\xi^\prime}^\beta \mathbf a(x,\,\xi^\prime,\,\tau)|\leq \langle \xi^\prime,\,\tau \rangle^{m-|\beta|+1} e^{-\delta^\prime \langle \xi^\prime,\,\tau \rangle |x_3-y_3|},\\
&& |(x_3-y_3) D_{x^\prime}^\alpha D_{\xi^\prime}^\beta \mathbf a(x,\,\xi^\prime,\,\tau)|\leq \langle \xi^\prime,\,\tau \rangle^{m-|\beta|-1} e^{-\delta^\prime \langle \xi^\prime,\,\tau \rangle |x_3-y_3|},\\
&& |D_{x^\prime} D_{x^\prime}^\alpha D_{\xi^\prime}^\beta \mathbf a(x,\,\xi^\prime,\,\tau)|\leq \langle \xi^\prime,\,\tau \rangle^{m-|\beta|} e^{-\delta^\prime \langle \xi^\prime,\,\tau \rangle |x_3-y_3|}
\end{eqnarray*}
for any $(x',\,\xi',\,\tau)\in \mathbb R^2_{x'}\times \mathbb R^2_{\xi'} \times \{\Sigma(\xi^\prime)\cap\{|\tau|>1\}\}$ with some constant $\delta^\prime>0$.
\par
If $\mathbf a(x^\prime,\,\xi^\prime,\,\tau)=\mathbf a(x^\prime,\,\xi^\prime,\,\tau;\,x_3,\,y,\,s)$ depends on the parameters $x_3\in\mathbb R_-,\,y\in{\mathbb R}_-^3,\,s\in(0,T)$, the above estimates have to be held uniformly with respect to them.
\end{defn}

We introduce the notations $\mathbf a_\ell,\,\mathbf b_\ell,\,\mathbf d_\ell,\,\mathbf e_\ell$ and $\mathcal A_\ell,\,\mathcal B_\ell,\,\mathcal D_\ell,\,\mathcal E_\ell$ as follows.

\medskip
$\mathbf a_\ell=\sum _{L=0}^{\infty}\mathbf a_{\ell,-1-L}$ with $e^{iy^\prime\cdot \xi^\prime}\mathbf a_{\ell,-1-L}\in S(-1-L)$ is the amplitude of the pseudo-differential operator $G_\ell^+$ defined by
\begin{equation}\label{p_G}
(G_\ell^+ \varphi) (x) = (2\pi)^{-2} \int_{\mathbb R^2}\int_{\mathbb R^2} e^{ix^\prime\cdot\xi^\prime} \mathbf a_\ell(x,\,\xi^\prime,\,\tau;\,y,\,s) \varphi(y^\prime)\,dy^\prime d\xi^\prime,\qquad \varphi(y^\prime)\in C_0^\infty({\mathbb R}^2).
\end{equation}
It holds that
\begin{equation*}
0=(\mathcal P G_\ell^+\varphi)(x) = (2\pi)^{-2}\int_{\mathbb R^2}\int_{\mathbb R^2} e^{ix'\cdot \xi'}\mathcal A_\ell (x,\,\xi',\,\tau;\,y,\,s)\varphi(y')\,dy^\prime d\xi'
\end{equation*}
with
$$\mathcal A_\ell (x,\,\xi',\,\tau;\,y,\,s)=\sum_{\alpha}(\alpha!)^{-1}\partial_{\xi'}^{\alpha}\mathcal P(x,\,\xi',\,\tau)D_{x'}^{\alpha}\mathbf a_\ell(x,\,\xi',\,\tau;\,y,\,s).$$
The Schwartz kernel $G_\ell^+(x,\,y)$ with $y=(y',\,y_3)\in\mathbb R_-^3,\,\, y_3<x_3\le0$ and the suppressed parameter $\tau$ is given by
\begin{equation}\label{kernel}
G_\ell^+(x,\,y)=(2\pi)^{-2} \int_{\mathbb R^2} e^{ix^\prime\cdot \xi^\prime} \mathbf a_\ell(x,\,\xi^\prime,\,\tau;\,y,\,s)\,d\xi^\prime.
\end{equation}
We will use a truncated sum for $\mathbf a_\ell$. Note that \eqref{kernel} differs from the usual definition of Schwartz kernel, but we shall see later that $\mathbf a_\ell(x,\,\xi^\prime,\,\tau;\,y,\,s)$ has the factor $e^{-iy^\prime\cdot\xi^\prime}$. Actually, each $e^{iy^\prime\cdot \xi^\prime}\mathbf a_{\ell,-1-L}\in S(-1-L)$ has a further property. That is, it can be analytically extended with respect to $(\xi^\prime,\,\tau)$ with $\tau=i\eta$ to $L^2_\mu$ with some estimates (see Section \ref{estimate}).

We arrange $\mathcal A_\ell (x,\,\xi',\,\tau;\,y,\,s)$ in terms of the order as follows:
\begin{equation*}
\mathcal A_\ell =\sum_{l=0}^{\infty}\mathcal A _{\ell,1-l},
\end{equation*}
where $\mathrm{ord}\,\mathcal A_{\ell,1-l}=1-l$, and $\mathcal A _{\ell,1-l}$ can be explicitly expressed; see, for example,
\begin{align*}
\mathcal A _{\ell,1} & = p_{2,0}^{(0)}\mathbf a_{\ell,-1},\\
\mathcal A _{\ell,0} & = \sum_{j+k+|\alpha|=1}(x_3-y_3)^j
p_{2,j}^{(\alpha)}D_{x'}^{\alpha}\mathbf a_{\ell,-1-k}+p_{1,0}^{(0)}\mathbf a_{\ell,-1}\\
& = p_{2,0}^{(0)}\mathbf a_{\ell,-2}+\sum_{j+|\alpha|=1}(x_3-y_3)^j
p_{2,j}^{(\alpha)}D_{x'}^{\alpha}\mathbf a_{\ell,-1}+p_{1,0}^{(0)}\mathbf a_{\ell,-1}
\end{align*}
with $p_{l,j}^{(\alpha)}=\partial_{x_3}^j\partial_{\xi'}^{\alpha}p_l(x',\,y_3,\,\xi',\,\tau)$, $l=1,\,2,\,j\in
\mathbb{Z}_{+},\, \alpha\in \mathbb{Z}_{+}^2$. We require $\{\mathbf a_{\ell,-1-L}\}_{L=0}^{\infty}$ to satisfy $\mathcal A_{\ell,1-l}=0$ for $l=0,\,1,\,\cdots$.

\medskip
Similarly, we define the amplitudes $\mathbf b_\ell,\,\mathbf d_\ell$ and $\mathbf e_\ell$ of $G_\ell^-,\,H_\ell^+$ and $H_\ell^-$. Also, corresponding to $\mathcal A_\ell $, we define the notations $\mathcal B_\ell,\,\mathcal D_\ell $ and $\mathcal E_\ell$ for $\mathbf b_\ell,\,\mathbf d_\ell$ and $\mathbf e_\ell$, respectively. Then we have the following representations of the amplitudes.

\begin{theorem}\label{param} Let $M^1=M|_{x_3=y_3}$ and $M^0=M|_{x_3=0}$. Apply these notations even for their components, for example, $m^1_{33}=m_{33}|_{x_3=y_3},\,m^0_{33}=m_{33}|_{x_3=0}$. Define
\begin{equation}\label{roots}
\begin{cases}
\lambda_{\pm} := (m_{33}^1)^{-1}\left[-i\displaystyle\sum_{j=1}^{2}m_{3j}^1\xi_j \pm \sqrt{m_{33}^1 (
\sum_{j,l=1}^{2}m_{jl}^1\xi_j\xi_l + \tau )-(\sum_{j=1}^2m_{3j}^1\xi_j)^2}\right], \\
\quad \\
\mu_{\pm} := (m_{33}^1)^{-1}\left[-i\displaystyle\sum_{j=1}^{2}m_{3j}^1\xi_j \pm \sqrt{m_{33}^1 (
\sum_{j,l=1}^{2}m_{jl}^1\xi_j\xi_l + \frac{\tau}{k})-(\sum_{j=1}^2m_{3j}^1\xi_j)^2}\right],
\end{cases}
\end{equation}
where the real parts of the square roots in $\lambda_{\pm}$ and $\mu_{\pm}$ are positive. Then we have
\begin{equation*}
\mathbf a_\ell=\sum_{L=0}^{\infty}\mathbf a_{\ell,-1-L},\quad \mathbf b_\ell=\sum_{L=0}^{\infty}\mathbf b_{\ell,-1-L},\quad
\mathbf d_\ell=\sum_{L=0}^{\infty}\mathbf d_{\ell,-1-L},\quad \mathbf e_\ell=\sum_{L=0}^{\infty}\mathbf e_{\ell,-1-L}
\end{equation*}
for $\ell=1,\,2$,
where
\begin{eqnarray}\label{amplitude1}
\mathbf a_{\ell,-L}&=&(2-\ell) \sum_{l=0}^{2L-2} f_{l,1}^{L}(x_3-y_3)^l\exp\left(\lambda_{+} x_3
- \lambda_{-} y_3-\tau s-iy'\cdot\xi'\right) \nonumber \\
&& +(\ell-1)\sum_{l=0}^{2L-2} f_{l,2}^{L}(x_3-y_3)^l
\exp\left(\lambda_{+} x_3 - \mu_{-} y_3-\tau s-iy'\cdot\xi'\right) \nonumber \\
&& +(2-\ell)\sum_{l=0}^{2L-2} f_{l,3}^{L}(x_3-y_3)^l \exp\left(\lambda_{-} x_3 - \lambda_{-} y_3-\tau s-iy'\cdot\xi'\right),
\end{eqnarray}
\begin{eqnarray}\label{amplitude12}
\mathbf b_{\ell,-L}&=&(2-\ell)\sum_{l=0}^{2L-2} f_{l,1}^{L}(x_3-y_3)^l\exp\left(\lambda_{+} x_3 -
\lambda_{-} y_3-\tau s-iy'\cdot\xi'\right) \nonumber\\
&& +(\ell-1)\sum_{l=0}^{2L-2} f_{l,2}^{L}(x_3-y_3)^l
\exp\left(\lambda_{+} x_3 - \mu_{-} y_3-\tau s-iy'\cdot\xi'\right) \nonumber\\
&&+(2-\ell)\sum_{l=0}^{2L-2}
f_{l,4}^{L}(x_3-y_3)^l \exp\left(\lambda_{+} x_3 - \lambda_{+} y_3-\tau s-iy'\cdot\xi'\right),
\end{eqnarray}
\begin{eqnarray}\label{amplitude13}
\mathbf d_{\ell,-L}&=&(2-\ell)\sum_{l=0}^{2L-2} f_{l,5}^{L}(x_3-y_3)^l \exp\left(\mu_{+} x_3 -
\lambda_{-} y_3-\tau s-iy'\cdot\xi'\right) \nonumber\\
& & + (\ell-1)\sum_{l=0}^{2L-2} f_{l,6}^{L}(x_3-y_3)^l
\exp\left(\mu_{+} x_3 - \mu_{-} y_3-\tau s-iy'\cdot\xi'\right) \nonumber \\
&& +(\ell-1) \sum_{l=0}^{2L-2} f_{l,7}^{L}(x_3-y_3)^l \exp\left(\mu_{-} x_3 - \mu_{-} y_3-\tau s-iy'\cdot\xi'\right),
\end{eqnarray}
and
\begin{eqnarray}\label{amplitude2}
\mathbf e_{\ell,-L}&=& (2-\ell)\sum_{l=0}^{2L-2} f_{l,5}^{L}(x_3-y_3)^l \exp\left(\mu_{+} x_3 -
\lambda_{-} y_3-\tau s-iy'\cdot\xi'\right)\nonumber \\
&& + (\ell-1)\sum_{l=0}^{2L-2} f_{l,6}^{L}(x_3-y_3)^l \exp\left(\mu_{+} x_3 - \mu_{-} y_3-\tau s-iy'\cdot\xi'\right) \nonumber \\
&& + (\ell-1)\sum_{l=0}^{2L-2} f_{l,8}^{L}(x_3-y_3)^l \exp\left(\mu_{+} x_3 -
\mu_{+} y_3-\tau s-iy'\cdot\xi'\right)
\end{eqnarray} with $\mathrm{ord}\, f_{l,j}^L = l-L$ for $j=1,\cdots,8$. In the above formulae, all $f_{l,j}^L$'s can be explicitly given in the following proof. Actually, $f_{l,j}^L$'s for $L=1$ are given in \eqref{Ord-1Amplitude} with $A_1,\,A_2,\,B_1,\,B_2$ defined in \eqref{coefficients}. For general $L\geq2$,  $f_{l,j}^L$'s are determined by $F_{l,j}$ in \eqref{coefficients-L1}--\eqref{coefficients-L2} and $A_5,\,A_7,\,B_5,\,B_7$ defined in \eqref{second_1}, \eqref{second_3}.
\end{theorem}

{\bf Proof.} We prove the result by induction on $L$. At first, let us find $\mathbf a_{\ell,-1},\,
\mathbf b_{\ell,-1},\,\mathbf d_{\ell,-1}$ and $\mathbf e_{\ell,-1}$. It implies from $\mathcal A_{\ell,1}=\mathcal B_{\ell,1}=\mathcal D_{\ell,1}=\mathcal E_{\ell,1}=0$ that
\begin{equation}\label{first}
p_{2,0}^{(0)} \mathbf a_{\ell,-1}=p_{2,0}^{(0)}\mathbf b_{\ell,-1}=q_{2,0}^{(0)}\mathbf d_{\ell,-1}=q_{2,0}^{(0)}\mathbf e_{\ell,-1}=0,
\end{equation}
where
\begin{eqnarray*}
p_{2,0}^{(0)} & = & m_{33}^1\partial_{x_3}^2+2i\sum_{j=1}^{2}m_{3j}^1\xi_j\partial_{x_3}- \left(\sum_{j,l=1}^{2}m_{jl}^1\xi_j\xi_l+\tau \right),\\
q_{2,0}^{(0)} & = & km_{33}^1\partial_{x_3}^2+2ik\sum_{j=1}^{2}m_{3j}^1\xi_j\partial_{x_3}-\left(k\sum_{j,l=1}^{2}m_{jl}^1\xi_j\xi_l + \tau \right).
\end{eqnarray*}
The solutions to the above ordinary differential equations \eqref{first} can be expressed as
\begin{eqnarray*}
\mathbf a_{\ell,-1} & = & C_1 \exp\left(\lambda_{+} x_3\right) + C_2 \exp\left(\lambda_{-} x_3\right),\\
\mathbf b_{\ell,-1} & = & C_3 \exp\left(\lambda_{+} x_3\right),\\
\mathbf d_{\ell,-1} & = & C_4 \exp\left(\mu_{+} x_3\right) + C_5 \exp\left(\mu_{-} x_3\right),\\
\mathbf e_{\ell,-1} & = & C_6 \exp\left(\mu_{+} x_3\right).
\end{eqnarray*}
Notice here that we took $\mathbf b_{\ell,-1}$ and $\mathbf e_{\ell,-1}$ satisfying $\displaystyle \lim_{x_3 \rightarrow -\infty}\mathbf b_{\ell,-1}=\lim_{x_3 \rightarrow -\infty}\mathbf e_{\ell,-1}=0$. From the transmission conditions and boundary conditions, we have
\begin{eqnarray*}
&& \mathbf a_{\ell,-1}-\mathbf b_{\ell,-1}=0,\quad ie_3\cdot
M^1\begin{pmatrix}\xi'\\D_3\end{pmatrix}(\mathbf a_{\ell,-1}-\mathbf b_{\ell,-1})=-(2-\ell)\mathcal J^{-1}(y)\exp\left(-\tau s-iy'\cdot\xi'\right)\qquad
{\rm on}\; x_3=y_3,\\
&& \mathbf d_{\ell,-1}-\mathbf e_{\ell,-1}=0,\quad ike_3\cdot
M^1\begin{pmatrix}\xi'\\D_3\end{pmatrix}(\mathbf d_{\ell,-1}-\mathbf e_{\ell,-1})=-(\ell -1)\mathcal J^{-1}(y)\exp\left(-\tau s-iy'\cdot\xi'\right)\qquad
{\rm on} \; x_3=y_3,\\
&& \mathbf a_{\ell,-1}-\mathbf d_{\ell,-1}=0,\quad ie_3\cdot
M^0\begin{pmatrix}\xi'\\D_3\end{pmatrix}\mathbf a_{\ell,-1}=ike_3\cdot
M^0\begin{pmatrix}\xi'\\D_3\end{pmatrix}\mathbf d_{\ell,-1}\qquad {\rm
on}\; x_3=0.
\end{eqnarray*}
From these conditions, we can easily derive the following system of equations
for constants $C_j\,(1\leq j\leq 6)$:
\begin{eqnarray}\label{A1_1}
&&C_1 \exp\left(\lambda_{+} y_3\right) + C_2 \exp\left(\lambda_{-} y_3\right) - C_3
\exp\left(\lambda_{+} y_3\right)=0,
\end{eqnarray}
\begin{eqnarray}\label{A1_2}
&&\lambda_{+} C_1 \exp\left(\lambda_{+} y_3\right) + \lambda_{-} C_2 \exp\left(\lambda_{-}
y_3\right) - \lambda_{+} C_3 \exp\left(\lambda_{+} y_3\right) \nonumber\\
&=& -(2-\ell)(m_{33}^1)^{-1}\mathcal J^{-1}(y)\exp\left(-\tau s-iy'\cdot\xi'\right),
\end{eqnarray}
\begin{eqnarray}\label{A1_3}
C_4 \exp\left(\mu_{+} y_3\right) + C_5 \exp\left(\mu_{-} y_3\right) - C_6 \exp\left(\mu_{+} y_3\right)=0,
\end{eqnarray}
\begin{eqnarray}\label{A1_4}
&&\mu_{+} C_4 \exp\left(\mu_{+} y_3\right) + \mu_{-} C_5 \exp\left(\mu_{-} y_3\right) - \mu_{+}
C_6 \exp\left(\mu_{+} y_3\right) \nonumber \\
& = & -(\ell-1)(km_{33}^1)^{-1}\mathcal J^{-1}(y)\exp\left(-\tau s-iy'\cdot\xi'\right),
\end{eqnarray}
\begin{eqnarray}\label{A1_5}
&&C_1 + C_2 - C_4 -C_5 =0,
\end{eqnarray}
\begin{eqnarray}\label{A1_6}
&&(i\displaystyle\sum_{j=1}^2m_{3j}^0\xi_j + \lambda_{+} m_{33}^0) C_1 +
(i\sum_{j=1}^2m_{3j}^0\xi_j + \lambda_{-} m_{33}^0) C_2 \nonumber\\
&=& k(i\displaystyle\sum_{j=1}^2m_{3j}^0\xi_j + \mu_{+} m_{33}^0) C_4 +
k(i\sum_{j=1}^2m_{3j}^0\xi_j + \mu_{-} m_{33}^0)C_5.
\end{eqnarray}
By \eqref{A1_1}--\eqref{A1_4}, we obtain
\begin{eqnarray*}
C_2 &=& \displaystyle\frac{(2-\ell)\mathcal J^{-1}(y)\exp\left(-\tau s-iy'\cdot\xi'\right)}{m_{33}^1(\lambda_{+}-\lambda_{-})}\exp\left(-\lambda_{-}y_3\right),\\
&&\\
C_5 &=& \displaystyle\frac{(\ell-1)\mathcal J^{-1}(y)\exp\left(-\tau s -iy'\cdot\xi'\right)}{km_{33}^1(\mu_{+}-\mu_{-})}\exp\left(-\mu_{-}y_3\right).
\end{eqnarray*}
Define
\begin{equation}\label{coefficients}
\left\{
\begin{array}{lll}
A_1&:=&\displaystyle\frac{\mathcal J^{-1}(y)}{m_{33}^1(\lambda_{+}-\lambda_{-})},\\
&&\\
B_1&:=&\displaystyle\frac{\mathcal J^{-1}(y)}{km_{33}^1(\mu_{+}-\mu_{-})},\\
&&\\
A_2&:=&\Big\{k(i\displaystyle\sum_{j=1}^2m_{3j}^0\xi_j+\mu_{+}m_{33}^0)-(i\sum_{j=1}^2m_{3j}^0\xi_j + \lambda_{+} m_{33}^0)\Big\}^{-1}\\
&&\quad \times (\lambda_{-}-\lambda_{+})m_{33}^0 A_1,\\
&&\\
B_2&:=&\Big\{k(i\displaystyle\sum_{j=1}^2m_{3j}^0\xi_j+\mu_{+}m_{33}^0)-(i\sum_{j=1}^2m_{3j}^0\xi_j + \lambda_{+} m_{33}^0)\Big\}^{-1}\\
&&\quad \{(i\displaystyle\sum_{j=1}^2m_{3j}^0\xi_j + \lambda_{+} m_{33}^0)-k(i\sum_{j=1}^2m_{3j}^0\xi_j + \mu_{-}m_{33}^0)\} B_1.
\end{array}
\right.
\end{equation}
Note that $\mathrm{ord}\, A_j=\mathrm{ord}\, B_j=-1$ for $j=1,\,2$. Then we derive from \eqref{A1_5} and
\eqref{A1_6} that
\begin{equation*}
C_4 =(2-\ell) A_2 \exp\left(-\lambda_{-} y_3-\tau s-iy'\cdot\xi'\right) + (\ell-1)B_2 \exp\left(-\mu_{-} y_3-\tau s-iy'\cdot\xi'\right),
\end{equation*}
and therefore
\begin{eqnarray*}
C_1 & = & C_4 - (2-\ell)A_1 \exp\left(-\lambda_{-} y_3-\tau s-iy'\cdot\xi'\right) + (\ell-1)B_1 \exp\left(-\mu_{-} y_3-\tau s-iy'\cdot\xi'\right)\\
& = & (2-\ell) (-A_1 + A_2) \exp\left(-\lambda_{-} y_3-\tau s-iy'\cdot\xi'\right) + (\ell-1)(B_1 + B_2) \exp\left(-\mu_{-}y_3-\tau s-iy'\cdot\xi'\right).
\end{eqnarray*}
Thus, we have
\begin{eqnarray*}
C_3 & = & (2-\ell) (-A_1 + A_2) \exp\left(-\lambda_{-} y_3-\tau s-iy'\cdot\xi'\right) + (\ell-1)(B_1 + B_2) \exp\left(-\mu_{-}
y_3-\tau s-iy'\cdot\xi'\right) \\
& & + (2-\ell)A_1 \exp\left(-\lambda_{+} y_3-\tau s-iy'\cdot\xi'\right),\\
C_6 & = & (2-\ell)A_2 \exp\left(-\lambda_{-} y_3-\tau s-iy'\cdot\xi'\right) + (\ell-1)B_2 \exp\left(-\mu_{-} y_3-\tau s-iy'\cdot\xi'\right)\\
& &  + (\ell-1)B_1\exp\left(-\mu_{+} y_3-\tau s-iy'\cdot\xi'\right).
\end{eqnarray*}
So we finally obtain
\begin{equation}\label{Ord-1Amplitude}
\left\{
\begin{array}{lll}
\mathbf a_{\ell,-1} & = & (2-\ell)(-A_1 + A_2) \exp\left(\lambda_{+} x_3 - \lambda_{-} y_3-\tau s-iy'\cdot\xi'\right)\\
& &  + (\ell-1)(B_1 + B_2)\exp\left(\lambda_{+} x_3 - \mu_{-} y_3-\tau s-iy'\cdot\xi'\right) \\
& &  + (2-\ell) A_1 \exp\left(\lambda_{-} x_3 - \lambda_{-} y_3-\tau s-iy'\cdot\xi'\right),\\
&&\\
\mathbf b_{\ell,-1} & = & (2-\ell)(-A_1 + A_2) \exp\left(\lambda_{+} x_3 - \lambda_{-} y_3-\tau s-iy'\cdot\xi'\right)  \\
& & + (\ell-1)(B_1 + B_2) \exp\left(\lambda_{+} x_3 - \mu_{-} y_3-\tau s-iy'\cdot\xi'\right) \\
& &  + (2-\ell)A_1 \exp\left(\lambda_{+} x_3 - \lambda_{+} y_3-\tau s-iy'\cdot\xi'\right),\\
&&\\
\mathbf d_{\ell,-1} & = & (2-\ell) A_2 \exp\left(\mu_{+} x_3 - \lambda_{-} y_3-\tau s-iy'\cdot\xi'\right) \\
& & + (\ell-1) B_2 \exp\left(\mu_{+}
x_3 - \mu_{-} y_3-\tau s-iy'\cdot\xi'\right)\\
& & + (\ell-1)B_1 \exp\left(\mu_{-} x_3 - \mu_{-} y_3-\tau s-iy'\cdot\xi'\right),\\
&&\\
\mathbf e_{\ell,-1} & = &(2-\ell) A_2 \exp\left(\mu_{+} x_3 - \lambda_{-} y_3-\tau s-iy'\cdot\xi'\right) \\
& & + (\ell-1)B_2 \exp\left(\mu_{+}
x_3 - \mu_{-} y_3-\tau s-iy'\cdot\xi'\right) \\
& & + (\ell-1)B_1 \exp\left(\mu_{+} x_3 - \mu_{+} y_3-\tau s-iy'\cdot\xi'\right).
\end{array}
\right.
\end{equation}
These show that \eqref{amplitude1}--\eqref{amplitude2} are true for $L=1$.

\medskip
We would like to give an important remark here.

\begin{remark}\label{symbols in Laplace domain}${}$

\noindent{\rm{(i)}} Ignoring the transmission conditions at $x_3=y_3$, we have the Lopatinskii matrix for $C_1$ and $C_4$ whose determinant is non-zero, by setting $C_2=C_5=0$ in \eqref{A1_5} and \eqref{A1_6}. In order to clarify the principal part of parametrices $G_\ell=G_\ell^\sharp$ and $H_\ell=H_\ell^\sharp$, we have included the transmission conditions at $x_3=y_3$ for the fundamental solutions of the partial differential operators $\mathcal{P}$ and $\mathcal{Q}$, respectively.

\medskip
\noindent{\rm{(ii)}} $A_1 \exp\left(\lambda_{-} x_3 - \lambda_{-} y_3-\tau s-iy'\cdot\xi'\right)$ in $(\mathbf a_{\ell,-1},\,\mathbf b_{\ell,-1})$ and $B_1 \exp\left(\mu_{+} x_3 - \mu_{+} y_3-\tau s-iy'\cdot\xi'\right)$ in $(\mathbf d_{\ell,-1},\,\mathbf e_{\ell,-1})$ give the principal parts of the fundamental solutions of $\mathcal{P}$ and $\mathcal{Q}$, respectively.
\end{remark}

Next, let us prove for the case $L=2$. From $\mathcal A_{\ell,0}=\mathcal B_{\ell,0}=\mathcal D_{\ell,0}=\mathcal E_{\ell,0}=0$, we know that
\begin{eqnarray*}
&&p_{2,0}^{(0)}\mathbf a_{\ell,-2}=-\sum_{j+|\alpha|=1}(x_3-y_3)^j
p_{2,j}^{(\alpha)}D_{x'}^\alpha \mathbf a_{\ell,-1}-p_{1,0}^{(0)}\mathbf a_{\ell,-1}=:\Theta_{\mathbf a_\ell},\\
&&p_{2,0}^{(0)}\mathbf b_{\ell,-2}=-\sum_{j+|\alpha|=1}(x_3-y_3)^j
p_{2,j}^{(\alpha)}D_{x'}^\alpha \mathbf b_{\ell,-1}-p_{1,0}^{(0)}\mathbf b_{\ell,-1}=:\Theta_{\mathbf b_\ell},\\
&&q_{2,0}^{(0)}\mathbf d_{\ell,-2}=-\sum_{j+|\alpha|=1}(x_3-y_3)^j
q_{2,j}^{(\alpha)}D_{x'}^\alpha \mathbf d_{\ell,-1}-q_{1,0}^{(0)}\mathbf d_{\ell,-1}=:\Theta_{\mathbf d_\ell},\\
&&q_{2,0}^{(0)}\mathbf e_{\ell,-2}=-\sum_{j+|\alpha|=1}(x_3-y_3)^j
q_{2,j}^{(\alpha)}D_{x'}^\alpha \mathbf e_{\ell,-1}-q_{1,0}^{(0)}\mathbf e_{\ell,-1}=:\Theta_{\mathbf e_\ell}.
\end{eqnarray*}
Note that
\begin{equation*}
\Theta_{\mathbf a_\ell}  = -\Big[(x_3-y_3)\partial_{x_3} p_2|_{x_3=y_3} - i\sum_{j=1}^2 \partial_{\xi_j}p_2|_{x_3=y_3} \partial_{x_j} + p_1|_{x_3=y_3} \Big] \mathbf a_{\ell,-1},
\end{equation*}
and $\Theta_{\mathbf b_\ell},\,\Theta_{\mathbf d_\ell},\,\Theta_{\mathbf e_\ell}$ can be expressed analogously.
From the forms of $\mathbf a_{\ell,-1},\,\mathbf b_{\ell,-1},\,\mathbf d_{\ell,-1}$ and $\mathbf e_{\ell,-1}$, we have
\begin{eqnarray*}
\Theta_{\mathbf a_\ell} & = & (2-\ell) \sum_{l=0}^{1} E_{l,1}(x_3-y_3)^l \exp\left(\lambda_{+} x_3 -
\lambda_{-} y_3 - \tau s - i y^\prime \cdot \xi^\prime\right)\\
&& + (\ell-1)\sum_{l=0}^{1} E_{l,2}(x_3-y_3)^l \exp\left(\lambda_{+}
x_3  -\mu_{-} y_3 - \tau s - i y^\prime \cdot \xi^\prime\right)\\
&& + (2-\ell)\sum_{l=0}^{1}E_{l,3}(x_3-y_3)^l \exp\left(\lambda_{-} x_3 - \lambda_{-} y_3 - \tau s - i y^\prime \cdot \xi^\prime\right),
\end{eqnarray*}
\begin{eqnarray*}
\Theta_{\mathbf b_\ell} & = & (2-\ell) \sum_{l=0}^{1} E_{l,1}(x_3-y_3)^l \exp\left(\lambda_{+} x_3 -
\lambda_{-} y_3 - \tau s - i y^\prime \cdot \xi^\prime\right) \\
&& + (\ell-1) \sum_{l=0}^{1} E_{l,2}(x_3-y_3)^l \exp\left(\lambda_{+}
x_3 - \mu_{-} y_3 - \tau s - i y^\prime \cdot \xi^\prime\right)\\
&& + (2-\ell) \sum_{l=0}^{1}
E_{l,4}(x_3-y_3)^l \exp\left(\lambda_{+} x_3 - \lambda_{+} y_3 - \tau s - i y^\prime \cdot \xi^\prime\right),
\end{eqnarray*}
\begin{eqnarray*}
\Theta_{\mathbf d_\ell} & = & (2-\ell) \sum_{l=0}^{1} E_{l,5}(x_3-y_3)^l \exp\left(\mu_{+} x_3 -
\lambda_{-} y_3 - \tau s - i y^\prime \cdot \xi^\prime\right)\\
&& + (\ell-1) \sum_{l=0}^{1} E_{l,6}(x_3-y_3)^l \exp\left(\mu_{+} x_3
- \mu_{-} y_3 - \tau s - i y^\prime \cdot \xi^\prime\right) \\
&& + (\ell-1) \sum_{l=0}^{1} E_{l,7}(x_3-y_3)^l \exp\left(\mu_{-} x_3 - \mu_{-} y_3 - \tau s - i y^\prime \cdot \xi^\prime\right),
\end{eqnarray*}
\begin{eqnarray*}
\Theta_{\mathbf e_\ell} & = & (2-\ell) \sum_{l=0}^{1} E_{l,5}(x_3-y_3)^l \exp\left(\mu_{+} x_3 -
\lambda_{-} y_3 - \tau s - i y^\prime \cdot \xi^\prime\right) \\
&& + (\ell-1) \sum_{l=0}^{1} E_{l,6}(x_3-y_3)^l \exp\left(\mu_{+} x_3
- \mu_{-} y_3 - \tau s - i y^\prime \cdot \xi^\prime\right) \\
&& + (\ell-1) \sum_{l=0}^{1} E_{l,8}(x_3-y_3)^l \exp\left(\mu_{+} x_3 - \mu_{+} y_3 - \tau s - i y^\prime \cdot \xi^\prime\right),
\end{eqnarray*}
where $E_{l,j}$ can be computed explicitly and $\mathrm{ord}\, E_{l,j}=l$ for $l=0,\,1$ and $j=1,\,\cdots,\,8$.

According to the form of $\Theta_{\mathbf a_\ell}$, $\mathbf a_{\ell,-2}$ can be expressed as
\begin{equation*}
\mathbf a_{\ell,-2}=\sum_{j=1}^3 \mathbf a_{\ell,-2}^j,\quad \mathbf a_{\ell,-2}^j= \alpha_j\sum_{l=0}^2
F_{l,j}(x_3-y_3)^l \exp\left(\beta_j x_3 - \delta_j y_3 - \tau s -iy^\prime\cdot\xi^\prime\right) \quad \mbox{for }
j=1,\,2,\,3
\end{equation*}
satisfying
\begin{equation*}
p_{2,0}^{(0)}\mathbf a_{\ell,-2}^j=\alpha_j\sum_{l=0}^1E_{l,j}(x_3-y_3)^l \exp\left(\beta_j x_3 - \delta_j y_3- \tau s -iy^\prime\cdot\xi^\prime\right),
\end{equation*}
where $\alpha_1=\alpha_3=2-\ell$, $\alpha_2=\ell -1$, $\beta_1=\beta_2=\lambda_{+}$, $\beta_3=\delta_1=\delta_3=\lambda_{-}$ and $\delta_2=\mu_{-}$.
From the above equation, we can express
$F_{l,j}$ in terms of $E_{l,j}$:
\begin{eqnarray}
F_{2,j} & = & \frac{E_{1,j}}{4\gamma_j(\beta_j
m_{33}^1+i\displaystyle\sum_{j=1}^2m_{3j}^1\xi_j)}, \label{coefficients_1}\\
F_{1,j} & = & \frac{E_{0,j}- 2\gamma_jm_{33}^1 F_{2,j}}{2\gamma_j(\beta_j
m_{33}^1+i\displaystyle\sum_{j=1}^2m_{3j}^1\xi_j)},
\label{coefficients_2}
\end{eqnarray}
where
\begin{equation*}
\gamma_j:=
\begin{cases}
1, & 1\leq j\leq 4,\\
k, & 5\leq j\leq 8,
\end{cases}
\end{equation*}
and $\mathrm{ord}\, F_{l,j}=l-2$ for $l=1,\,2$. We do the same calculations for $\mathbf b_{\ell,-2},\, \mathbf d_{\ell,-2}$ and $\mathbf e_{\ell,-2}$. Then we have
\begin{eqnarray}\label{L2_1}
\mathbf a_{\ell,-2}&=& (2-\ell) \Big\{\sum_{l=1}^{2} F_{l,1}(x_3-y_3)^l +
C_1\Big\}\exp\left(\lambda_{+} x_3 - \lambda_{-} y_3- \tau s -iy^\prime\cdot\xi^\prime\right) \nonumber\\
& & +(\ell-1)\Big\{\sum_{l=1}^{2}
F_{l,2}(x_3-y_3)^l + C_2\Big\}\exp\left(\lambda_{+} x_3 - \mu_{-} y_3- \tau s -iy^\prime\cdot\xi^\prime\right) \nonumber\\
& & +(2-\ell)\Big\{\sum_{l=1}^{2}
F_{l,3}(x_3-y_3)^l + C_3\Big\} \exp\left(\lambda_{-} x_3 - \lambda_{-} y_3- \tau s -iy^\prime\cdot\xi^\prime\right),
\end{eqnarray}
\begin{eqnarray}\label{L2_2}
\mathbf b_{\ell,-2}&=&(2-\ell)\Big\{\sum_{l=1}^{2} F_{l,1}(x_3-y_3)^l +
C_4\Big\}\exp\left(\lambda_{+} x_3 - \lambda_{-} y_3- \tau s -iy^\prime\cdot\xi^\prime\right)  \nonumber\\
&& +(\ell-1)\Big\{\sum_{l=1}^{2}
F_{l,2}(x_3-y_3)^l + C_5\Big\} \exp\left(\lambda_{+} x_3 - \mu_{-} y_3- \tau s -iy^\prime\cdot\xi^\prime\right) \nonumber\\
&& +(2-\ell)\Big\{\sum_{l=1}^{2}
F_{l,4}(x_3-y_3)^l + C_6\Big\} \exp\left(\lambda_{+} x_3 - \lambda_{+} y_3- \tau s -iy^\prime\cdot\xi^\prime\right),
\end{eqnarray}
\begin{eqnarray}\label{L2_3}
\mathbf d_{\ell,-2}&=&(2-\ell)\Big\{\sum_{l=1}^{2} F_{l,5}(x_3-y_3)^l +
C_7\Big\}\exp\left(\mu_{+} x_3 - \lambda_{-} y_3- \tau s -iy^\prime\cdot\xi^\prime\right)  \nonumber\\
&& + (\ell-1)\Big\{\sum_{l=1}^{2}
F_{l,6}(x_3-y_3)^l + C_8\Big\}\exp\left(\mu_{+} x_3 - \mu_{-} y_3- \tau s -iy^\prime\cdot\xi^\prime\right) \nonumber\\
& & + (\ell-1)\Big\{\sum_{l=1}^{2}
F_{l,7}(x_3-y_3)^l +C_9\Big\}\exp\left(\mu_{-} x_3 - \mu_{-} y_3- \tau s -iy^\prime\cdot\xi^\prime\right),
\end{eqnarray}
\begin{eqnarray}\label{L2_4}
\mathbf e_{\ell,-2}&=& (2-\ell)\Big\{\sum_{l=1}^{2} F_{l,5}(x_3-y_3)^l +
C_{10}\Big\}\exp\left(\mu_{+} x_3 - \lambda_{-} y_3- \tau s -iy^\prime\cdot\xi^\prime\right)  \nonumber\\
&& + (\ell-1)\Big\{\sum_{l=1}^{2}
F_{l,6}(x_3-y_3)^l + C_{11}\Big\}\exp\left(\mu_{+} x_3 - \mu_{-} y_3- \tau s -iy^\prime\cdot\xi^\prime\right) \nonumber\\
&& + (\ell-1)\Big\{\sum_{l=1}^{2} F_{l,8}(x_3-y_3)^l + C_{12}\Big\}\exp\left(\mu_{+} x_3 - \mu_{+} y_3- \tau s -iy^\prime\cdot\xi^\prime\right),
\end{eqnarray}
where $C_j\,(j=1,\,\cdots,\,12)$ are constants with respect to $x_3$. From the following transmission and boundary conditions:
\begin{eqnarray}
&&\mathbf a_{\ell,-2}-\mathbf b_{\ell,-2}=0,\quad ie_3\cdot M^1\begin{pmatrix}\xi'\\D_3\end{pmatrix}(\mathbf a_{\ell,-2}-\mathbf b_{\ell,-2})=0\quad {\rm on} \; x_3=y_3, \label{L2_5}\\
&&\mathbf d_{\ell,-2}-\mathbf e_{\ell,-2}=0,\quad ike_3\cdot M^1\begin{pmatrix}\xi'\\D_3\end{pmatrix}(\mathbf d_{\ell,-2}-\mathbf e_{\ell,-2})=0\quad {\rm on} \; x_3=y_3, \label{L2_6}\\
&&\mathbf a_{\ell,-2}-\mathbf d_{\ell,-2}=0,\quad ie_3\cdot M^0\begin{pmatrix}\xi'\\D_3\end{pmatrix}\mathbf a_{\ell,-2}=ike_3\cdot
M^0\begin{pmatrix}\xi'\\D_3\end{pmatrix}\mathbf d_{\ell,-2}\quad{\rm on}\; x_3=0, \label{L2_7}
\end{eqnarray}
we can derive a system of linear equations for $C_j \, (j=1,\,\cdots,\,12)$.

Looking at the structures of $\mathbf a_{\ell,-2},\,\mathbf b_{\ell,-2},\,\mathbf d_{\ell,-2}$ and $\mathbf e_{\ell,-2}$, we only have six substantial unknowns, namely, $C_3,\,C_9$, $(2-\ell)C_1 \exp\left(-\lambda_- y_3\right) + (\ell-1)C_2 \exp\left(-\mu_- y_3\right)$, $(2-\ell)C_4 \exp\left(-\lambda_- y_3\right)+ (\ell-1)C_5 \exp\left(-\mu_- y_3\right) + (2-\ell)C_6 \exp\left(-\lambda_+ y_3\right)$, $(2-\ell)C_7\exp\left(-\lambda_ - y_3\right) + (\ell-1)C_8 \exp\left(-\mu_- y_3\right)$ and $(2-\ell)C_{10} \exp\left(-\lambda_- y_3\right) + (\ell-1)C_{11}\exp\left(-\mu_- y_3\right) + (\ell-1)C_{12} \exp\left(-\mu_+ y_3\right)$. Hence,  $\mathbf a_{\ell,-2},\,\mathbf b_{\ell,-2},\,\mathbf d_{\ell,-2}$ and $\mathbf e_{\ell,-2}$ can be uniquely determined, and the details will be shown in the appendix.

\bigskip
Suppose that \eqref{amplitude1}--\eqref{amplitude2} are true for $L\geq 2$. We will show that they also hold for $L+1$. Note that
\begin{align}\label{A_L}
\mathcal A_{\ell,1-L} & = p_{2,0}^{(0)} \mathbf a_{\ell,-1-L} + \sum _{j+k+|\alpha|=L\atop k\leq
L-1}(\alpha!)^{-1}(j!)^{-1}(x_3-y_3)^jp_{2,j}^{(\alpha)}D_{x'}^{\alpha}\mathbf a_{\ell,-1-k} \nonumber\\
& \quad +
\sum_{j+k+|\alpha|=L-1}(\alpha!)^{-1}(j!)^{-1}(x_3-y_3)^jp_{1,j}^{(\alpha)}D_{x'}^{\alpha}\mathbf a_{\ell,-1-k} \nonumber\\
& \quad +
\sum_{j+k=L-2}(j!)^{-1}(x_3-y_3)^jp_{0,j}^{(0)}\mathbf a_{\ell,-1-k}.
\end{align}
Then $\mathcal A_{\ell,1-L}=0$ implies that
\begin{align*}
p_{2,0}^{(0)}\mathbf a_{\ell,-1-L} & =
-\sum_{k=0}^{L-1}\Big[\sum_{j+k+|\alpha|=L}(\alpha!)^{-1}(j!)^{-1}(x_3-y_3)^jp_{2,j}^{(\alpha)}D_{x'}^{\alpha}\mathbf a_{\ell,-1-k}\\
& \qquad \qquad \quad +
\sum_{j+k+|\alpha|=L-1}(\alpha!)^{-1}(j!)^{-1}(x_3-y_3)^jp_{1,j}^{(\alpha)}D_{x'}^{\alpha}\mathbf a_{\ell,-1-k}\Big]\\
& \quad -
\sum_{j+k=L-2}(j!)^{-1}(x_3-y_3)^jp_{0,j}^{(0)}\mathbf a_{\ell,-1-k}\\
&=:\Theta_{\mathbf a_\ell}^{1-L}.
\end{align*}
From the form of $\mathbf a_{\ell,-1-k}$, we know that $\Theta_{\mathbf a_\ell}^{1-L}$ is the sum of $(2-\ell)\exp\left(\lambda_{+} x_3 - \lambda_{-}
y_3 -\tau s -iy^\prime\cdot \xi^\prime\right)$, $(\ell-1)\exp\left(\lambda_{+} x_3 - \mu_{-} y_3 -\tau s -iy^\prime\cdot \xi^\prime\right)$ and $(2-\ell)\exp\left(\lambda_{-} x_3 - \lambda_{-} y_3 -\tau s -iy^\prime\cdot \xi^\prime\right)$ with $(2L-1)$-th degree's polynomials in $x_3-y_3$ as the coefficients of exponentials. That is,
\begin{eqnarray*}
\Theta_{\mathbf a_\ell}^{1-L} &=&  (2-\ell)\sum_{l=0}^{2L-1} E_{l,1}(x_3-y_3)^l \exp\left(\lambda_{+} x_3 - \lambda_{-} y_3 -\tau s -iy^\prime\cdot \xi^\prime\right) \\ &&+(\ell-1)\sum_{l=0}^{2L-1} E_{l,2}(x_3-y_3)^l \exp\left(\lambda_{+} x_3 - \mu_{-} y_3 -\tau s -iy^\prime\cdot \xi^\prime\right)\\
&&+(2-\ell)\sum_{l=0}^{2L-1} E_{l,3}(x_3-y_3)^l \exp\left(\lambda_{-} x_3 - \lambda_{-} y_3 -\tau s -iy^\prime\cdot \xi^\prime\right),
\end{eqnarray*}
where $\mathrm{ord}\,E_{l,j}=l+1-L$ for $j=1,\,2,\,3$. So we have
\begin{equation*}
\mathbf a_{\ell,-1-L}=\sum_{j=1}^3 \mathbf a_{\ell,-1-L}^j,\quad \mathbf a_{\ell,-1-L}^j=\alpha_j\sum_{l=0}^{2L}
F_{l,j}(x_3-y_3)^l \exp\left(\beta_j x_3 - \delta_j y_3 -\tau s -iy^\prime\cdot \xi^\prime\right),\quad j=1,\,2,\,3
\end{equation*}
satisfying
\begin{equation*}
p_{2,0}^{(0)}\mathbf a_{\ell,-1-L}^j=\alpha_j\sum_{l=0}^{2L-1}E_{l,j}(x_3-y_3)^l \exp\left(\beta_j x_3 - \delta_j y_3 -\tau s -iy^\prime\cdot \xi^\prime\right),
\end{equation*}
where $\alpha_1=\alpha_3=2-\ell$, $\alpha_2=\ell-1$, $\beta_1=\beta_2=\lambda_{+}, \,\beta_3=\delta_1=\delta_3=\lambda_{-}$ and $\delta_2=\mu_{-}$.
We can express $F_{l,j}$ in terms of $E_{l,j}$:
\begin{eqnarray}
F_{2L,j} & = & \frac{E_{2L-1,j}}{4L(\beta_j m_{33}^1+i\displaystyle\sum_{j=1}^2m_{3j}^1\xi_j)}, \label{coefficients-L1}\\
F_{l+1,j} & = & \frac{E_{l,j}- (l+2)(l+1)m_{33}^1 F_{l+2,j}}{2(l+1)(\beta_j m_{33}^1+i\displaystyle\sum_{j=1}^2m_{3j}^1\xi_j)}, \label{coefficients-L2}
\end{eqnarray}
where $\mathrm{ord}\, F_{l+1,j}=l-L$ for $l=0,\,1,\,\cdots,\,2L-1$. Therefore, we have
\begin{eqnarray*}
\mathbf a_{\ell,-1-L}&=& (2-\ell) \Big\{\sum_{l=1}^{2L} f_{l,1}^{L+1}(x_3-y_3)^l +
C_1\Big\}\exp\left(\lambda_{+} x_3 - \lambda_{-} y_3 -\tau s -iy^\prime\cdot \xi^\prime\right)\\
&& +(\ell-1)\Big\{\sum_{l=1}^{2L}f_{l,2}^{L+1}(x_3-y_3)^l + C_2\Big\}\exp\left(\lambda_{+} x_3 - \mu_{-} y_3 -\tau s -iy^\prime\cdot \xi^\prime\right)\\
& & +(2-\ell)\Big\{\sum_{l=1}^{2L} f_{l,3}^{L+1}(x_3-y_3)^l + C_3\Big\}
\exp\left(\lambda_{-} x_3 - \lambda_{-} y_3 -\tau s -iy^\prime\cdot \xi^\prime\right),
\end{eqnarray*}
where $f_{l,j}^{L+1}=F_{l,j}$ for $l=1,\,\cdots,\, 2L$ and $ j=1,\,2,\,3$. In the same way, we can get the similar expressions for $\mathbf b_{\ell,-1-L},\, \mathbf d_{\ell,-1-L}$ and $\mathbf e_{\ell,-1-L}$ as follows:
\begin{eqnarray*}
\mathbf b_{\ell,-1-L}&=&(2-\ell)\Big\{\sum_{l=1}^{2L} f_{l,1}^{L+1}(x_3-y_3)^l + C_4\Big\}\exp\left(\lambda_{+} x_3 - \lambda_{-} y_3 -\tau s -iy^\prime\cdot \xi^\prime\right)\\
&&+(\ell-1)\Big\{\sum_{l=1}^{2L} f_{l,2}^{L+1}(x_3-y_3)^l + C_5\Big\} \exp\left(\lambda_{+} x_3 - \mu_{-} y_3 -\tau s -iy^\prime\cdot \xi^\prime\right)\\
&& +(2-\ell)\Big\{\sum_{l=1}^{2L} f_{l,4}^{L+1}(x_3-y_3)^l + C_6\Big\} \exp\left(\lambda_{+} x_3 - \lambda_{+} y_3 -\tau s -iy^\prime\cdot \xi^\prime\right),
\end{eqnarray*}
\begin{eqnarray*}
\mathbf d_{\ell,-1-L}&=&(2-\ell)\Big\{\sum_{l=1}^{2L} f_{l,5}^{L+1}(x_3-y_3)^l + C_7\Big\}\exp\left(\mu_{+} x_3 - \lambda_{-} y_3 -\tau s -iy^\prime\cdot \xi^\prime\right) \\
&& + (\ell-1)\Big\{\sum_{l=1}^{2L} f_{l,6}^{L+1}(x_3-y_3)^l + C_8\Big\}\exp\left(\mu_{+} x_3 - \mu_{-} y_3 -\tau s -iy^\prime\cdot \xi^\prime\right)\\
&& + (\ell-1)\Big\{\sum_{l=1}^{2L} f_{l,7}^{L+1}(x_3-y_3)^l + C_9\Big\}\exp\left(\mu_{-} x_3 - \mu_{-} y_3 -\tau s -iy^\prime\cdot \xi^\prime\right),
\end{eqnarray*}
\begin{eqnarray*}
\mathbf e_{\ell,-1-L}&=&(2-\ell) \Big\{\sum_{l=1}^{2L} f_{l,5}^{L+1}(x_3-y_3)^l + C_{10}\Big\}\exp\left(\mu_{+} x_3 - \lambda_{-} y_3 -\tau s -iy^\prime\cdot \xi^\prime\right) \\
&& + (\ell-1)\Big\{\sum_{l=1}^{2L} f_{l,6}^{L+1}(x_3-y_3)^l + C_{11}\Big\}\exp\left(\mu_{+} x_3 - \mu_{-} y_3 -\tau s -iy^\prime\cdot \xi^\prime\right)\\
&& + (\ell-1)\Big\{\sum_{l=1}^{2L} f_{l,8}^{L+1}(x_3-y_3)^l + C_{12}\Big\}\exp\left(\mu_{+} x_3 - \mu_{+} y_3 -\tau s -iy^\prime\cdot \xi^\prime\right),
\end{eqnarray*}
where $f_{l,j}^{L+1}$ are determined from $E_{i,j}$ $(l=1,\,\cdots,\,2L,\; i=0,\,1,\,\cdots,\,2L-1,\; j=1,\,\cdots,\,8)$ in $\Theta_\mathbf m\,(\mathbf m=\mathbf b_\ell,\,\mathbf d_\ell,\,\mathbf e_\ell)$ like \eqref{coefficients-L1} and \eqref{coefficients-L2}. To determine the constants $C_j\,(j=1,\,\cdots,\,12)$, we use the transmission conditions on $x_3=y_3$ and the boundary conditions $x_3=0$, namely,
\begin{eqnarray*}
&&\mathbf a_{\ell,-1-L}-\mathbf b_{\ell,-1-L}=0,\quad ie_3\cdot M^1\begin{pmatrix}\xi'\\D_3\end{pmatrix}(\mathbf a_{\ell,-1-L}-\mathbf b_{\ell,-1-L})=0\qquad \mathrm{on} \; x_3=y_3,\\
&&\mathbf d_{\ell,-1-L}-\mathbf e_{\ell,-1-L}=0,\quad ike_3\cdot M^1\begin{pmatrix}\xi'\\D_3\end{pmatrix}(\mathbf d_{\ell,-1-L}-\mathbf e_{\ell,-1-L})=0\qquad \mathrm{on} \; x_3=y_3,\\
&&\mathbf a_{\ell,-1-L}-\mathbf d_{\ell,-1-L}=0,\quad ie_3\cdot M^0\begin{pmatrix}\xi'\\D_3\end{pmatrix}\mathbf a_{\ell,-1-L}=ike_3\cdot
M^0\begin{pmatrix}\xi'\\D_3\end{pmatrix}\mathbf d_{\ell,-1-L}\qquad \mathrm{on}\; x_3=0.
\end{eqnarray*}
These equations lead to a system of equations for $C_j\, (j=1,\,\cdots,\,12)$, which has only six substantial unknowns as for $L=2$. This system can be solved through the same process as we did for $L=2$ in Appendix 1. Then we have the following expressions:
\begin{eqnarray*}
\mathbf a_{\ell,-1-L} &= & (2-\ell)\Big\{\sum_{l=1}^{2L} f_{l,1}^{L+1}(x_3-y_3)^l +
(A_5+A_7)\Big\}\exp\left(\lambda_{+} x_3 - \lambda_{-} y_3 -\tau s - iy^\prime \cdot \xi^\prime\right)\\
&& + (\ell-1)\Big\{\sum_{l=1}^{2L}
f_{l,2}^{L+1}(x_3-y_3)^l + (B_5+B_7)\Big\}\exp\left(\lambda_{+} x_3 - \mu_{-} y_3 -\tau s - iy^\prime \cdot \xi^\prime\right)\\
&& +(2-\ell)\Big\{\sum_{l=1}^{2L} f_{l,3}^{L+1}(x_3-y_3)^l +
\frac{f_{1,3}^{L+1}-f_{1,4}^{L+1}}{\lambda_{+}-\lambda_{-}}\Big\}
\exp\left(\lambda_{-} x_3 - \lambda_{-} y_3 -\tau s - iy^\prime \cdot \xi^\prime\right),
\end{eqnarray*}
\begin{eqnarray*}
\mathbf b_{\ell,-1-L} &=& (2-\ell)\Big\{\sum_{l=1}^{2L} f_{l,1}^{L+1}(x_3-y_3)^l +
(A_5+A_7)\Big\}\exp\left(\lambda_{+} x_3 - \lambda_{-} y_3 -\tau s - iy^\prime \cdot \xi^\prime\right)\\
&& + (\ell-1)\Big\{\sum_{l=1}^{2L}
f_{l,2}^{L+1}(x_3-y_3)^l + (B_5+B_7)\Big\} \exp\left(\lambda_{+} x_3 - \mu_{-} y_3 -\tau s - iy^\prime \cdot \xi^\prime\right)\\
&& +(2-\ell)\Big\{\sum_{l=1}^{2L} f_{l,4}^{L+1}(x_3-y_3)^l +
\frac{f_{1,3}^{L+1}-f_{1,4}^{L+1}}{\lambda_{+}-\lambda_{-}}\Big\} \exp\left(\lambda_{+} x_3 - \lambda_{+} y_3 -\tau s - iy^\prime \cdot \xi^\prime\right),
\end{eqnarray*}
\begin{eqnarray*}
\mathbf d_{\ell,-1-L} & = & (2-\ell)\Big\{\sum_{l=1}^{2L} f_{l,5}^{L+1}(x_3-y_3)^l + A_7\Big\}\exp\left(\mu_{+} x_3 - \lambda_{-} y_3 -\tau s - iy^\prime \cdot \xi^\prime\right)\\
&& + (\ell-1)\Big\{\sum_{l=1}^{2L} f_{l,6}^{L+1}(x_3-y_3)^l + B_7\Big\}\exp\left(\mu_{+} x_3 - \mu_{-} y_3 -\tau s - iy^\prime \cdot \xi^\prime\right)\\
&& + (\ell-1)\Big\{\sum_{l=1}^{2L} f_{l,7}^{L+1}(x_3-y_3)^l +
\frac{f_{1,7}^{L+1} - f_{1,8}^{L+1}}{\mu_{+}-\mu_{-}}\Big\}\exp\left(\mu_{-} x_3 - \mu_{-} y_3 -\tau s - iy^\prime \cdot \xi^\prime\right),
\end{eqnarray*}
\begin{eqnarray*}
\mathbf e_{\ell,-1-L} & = &  (2-\ell)\Big\{\sum_{l=1}^{2L} f_{l,5}^{L+1}(x_3-y_3)^l + A_7\Big\}\exp\left(\mu_{+} x_3 - \lambda_{-} y_3 -\tau s - iy^\prime \cdot \xi^\prime\right)\\
&& + (\ell-1)\Big\{\sum_{l=1}^{2L} f_{l,6}^{L+1}(x_3-y_3)^l + B_7\Big\}\exp\left(\mu_{+} x_3 - \mu_{-} y_3 -\tau s - iy^\prime \cdot \xi^\prime\right)\\
&& + (\ell-1)\Big\{\sum_{l=1}^{2L} f_{l,8}^{L+1}(x_3-y_3)^l + \frac{f_{1,7}^{L+1} - f_{1,8}^{L+1}}{\mu_{+}-\mu_{-}}\Big\}\exp\left(\mu_{+} x_3 - \mu_{+}
y_3 -\tau s - iy^\prime \cdot \xi^\prime\right),
\end{eqnarray*}
where $\mathrm{ord}\,A_i=\mathrm{ord}\,B_i=-1-L$ for $i=5,\,7$ and $\mathrm{ord}\, f_{l,j}^{L+1}=l-1-L$ for $j=1,\,\cdots,\,8$. Thus, we have shown \eqref{amplitude1}--\eqref{amplitude2} for $L+1$. The proof of the theorem is complete. \hfill $\Box$

\bigskip
Once we have determined the amplitudes $\mathbf a_\ell,\,\mathbf b_\ell,\,\mathbf d_\ell,\,\mathbf e_\ell$, the corresponding pseudo-differential operators and their Schwartz kernels can be expressed in terms of \eqref{p_G} and \eqref{kernel}, respectively.

\section{Construction of a local parametrix for parabolic ITP}\label{estimate}
\setcounter{equation}{0}

In this section, we construct a local parametrix for the parabolic interior transmission problem \eqref{eq:ITP-2} by taking the inverse Laplace transform of $(G_\ell^\pm,\,H_\ell^\pm)$ given in the last section. The error estimates coming from the truncation of the amplitudes are derived. We only show how to handle $G_1^{+}$, since the arguments for the others are the same. To proceed, we need the following result given in \cite{Arima, Greiner1971}:
\begin{lemma}\label{Arima} For each $\rho\geq 0$, let $g(\xi',\,\eta,\,\rho)$ be a holomorphic function of $(\xi',\,\eta)$
in $L_{\mu}^2 \subset \mathbb{C}^2\times \mathbb{C}$ for
some $\mu>0$ with
\begin{equation*}
L_{\mu}^2=\left\{(\xi',\,\eta)\in
\mathbb{C}^2\times\mathbb{C}:\;\mathrm{Im}\,\eta < \mu(|\mathrm{Re}\, \eta|+|\mathrm{Re}\,
\xi'|^2)-\mu^{-1}|\mathrm{Im}\,\xi'|^2\right\}.\end{equation*}
Assume that
\begin{equation}\label{eq:estimate for g}
|g(\xi',\,\eta,\,\rho)|\leq C(|\xi'|+|\eta|^{1/2})^\kappa \exp[-c\rho
(|\xi'|+|\eta|^{1/2})]\end{equation} for $(\xi',\,\eta)\in
L_{\mu}^2$, $\kappa\leq 0$ and $\rho\geq 0$. Set
\begin{equation*}
\begin{cases}
G(x',\,t;\,\rho)=(2\pi)^{-2}\displaystyle\int_{{\mathbb R}^2}e^{ix'\cdot\xi'}\int_{-\infty-iq}^{\infty-iq}e^{it\eta}g(\xi',\,\eta,\,\rho)\,d\eta\,
d\xi', & \rho>0,\\
G(x',\,t;\,0)= \displaystyle\lim_{\rho\downarrow 0}G(x',\,t;\,\rho),&
\end{cases}
\end{equation*}
where $q=q(\xi^\prime)$ is an arbitrary positive number such that $(\xi^\prime,\,\eta)$ with $\eta=\gamma-iq\,\,(\gamma\in\mathbb R)$ is in $L^2_\mu$. Then we have
\begin{equation}\label{estimate_G}
|G(x',\,t;\,\rho)|\leq
c_1t^{-\frac{\kappa}{2}-2}\exp\Big[-c_2\frac{|x'|^2+{\rho}^2}{t}\Big],\qquad x^\prime\in\mathbb R^2,\,t>0,\,\rho\geq 0,
\end{equation}
where $c_1$ and $c_2$ are positive constants independent of $x^\prime,\,t,\,\rho$. Moreover, $g(\xi^\prime,\,\eta,\,\rho)$ can have parameters. Namely, let $g$ be analytic in $L^2_\mu$ with respect to $(\xi^\prime,\,\eta)$ by fixing these new parameters together with $\rho$ and satisfy the estimate \eqref{eq:estimate for g} uniformly with respect to the new parameters. Then we have the final estimate \eqref{estimate_G} uniformly with respect to the new parameters. \end{lemma}

\begin{remark}\label{Arima_R}
If $t$ is positive, then $G$ is smooth with respect to $x^\prime,\,t,\,\rho$. Further, if $g$ has parameters and depends smoothly on them without destroying the analyticity in $L^2_\mu$ and the estimate \eqref{eq:estimate for g} likewise stated in the above lemma, then $G$ is also smooth with respect to the parameters. These are well known and can be easily checked by looking at the proof of Lemma \ref{Arima}.
\end{remark}

Recall the Schwartz kernel \eqref{kernel} that
\begin{eqnarray}\label{w1}
G_1^+(x,\,y) & = &(2\pi)^{-2}\int_{\mathbb R^2}e^{i x'\cdot
\xi'}\mathbf a_1(x,\,\xi',\,\tau;\,y,\,s)\,d\xi' \nonumber\\
& = & (2\pi)^{-2}\int_{\mathbb R^2}e^{i x'\cdot
\xi'}\sum_{j=0}^{\infty}\mathbf a_{1,-1-j}(x,\,\xi',\,\tau;\,y,\,s)\,d\xi'.
\end{eqnarray}
Since we will truncate $\sum_{j=0}^{\infty}\mathbf a_{1,-1-j}$ at $j=N$ with large $N$, we do not care about the convergence of the formal sum.  $G_\ell^\pm$ and $H_\ell^\pm$ for $\ell=1,\,2$ with the amplitudes truncated at $j=N$ are denoted by $G_{\ell,N}^\pm$ and $H_{\ell,N}^\pm$, respectively. $\dot G_{\ell,N}^\pm(x,\,t;\,y,\,s)$ and $\dot H_{\ell,N}^\pm(x,\,t;\,y,\,s)$ are the inverse Laplace transforms of $G_{\ell,N}^\pm(x,\,\tau;\,y,\,s)$ and $H_{\ell,N}^\pm(x,\,\tau;\,y,\,s)$ with respect to $\tau$. Concerning the estimate of the inverse Laplace transform of $G_1^+(x,\,y) = G_1^+(x,\,y;\,\tau,\,s)$ with respect to $\tau$, we will only show the estimate for the leading term of \eqref{w1}, which gives the dominant part of the estimate. Let
\begin{equation*}
\dot G_1^+(x,\,t;\,y,\,s)=(2\pi i)^{-1}(2\pi)^{-2}\int_{\sigma-i\infty}^{\sigma+i\infty}\int_{\mathbb R^2}\exp\left(t\tau
+ ix'\cdot\xi'\right)\mathbf a_{1,-1}(x,\,\xi',\,\tau;\,y,\,s)\,d\xi'd\tau,
\end{equation*}
where
\begin{eqnarray*}
\mathbf a_{1,-1} & = & (-A_1 + A_2) \exp\left(\lambda_{+} x_3 - \lambda_{-} y_3-\tau s-iy'\cdot\xi'\right)\\
& & + A_1 \exp\left(\lambda_{-} x_3 - \lambda_{-} y_3-\tau s-iy'\cdot\xi'\right)\\
&=:& \mathbf a_{1,-1}^{(1)} + \mathbf a_{1,-1}^{(2)}.
\end{eqnarray*}

Define
\begin{equation}\label{dot G^+_1(29)}
\dot G^+_{1(j)}(x,\,t;\,y,\,s):=-(2\pi)^{-3}i\int_{\sigma-i\infty}^{\sigma+i\infty}\int_{\mathbb R^2}\exp\left(t\tau
+ ix'\cdot\xi'\right) \mathbf a_{1,-1}^{(j)}(x,\,\xi',\,\tau;\,y,\,s)\,d\xi'd\tau,\quad j=1,\,2.
\end{equation}
We only estimate $\dot G^+_{1(1)}$, since $\dot G^+_{1(2)}$ can be done analogously. By the change of variable $\tau=i\eta$, we have
\begin{eqnarray*}
\dot G^+_{1(1)}(x,\,t;\,y,\,s)
& = & -(2\pi)^{-3}i\int_{\sigma-i\infty}^{\sigma+i\infty}\int_{\mathbb R^2}\exp\left(t\tau
+ ix'\cdot\xi'\right) (-A_1+A_2) \\
& & \qquad\qquad\qquad\qquad\qquad\times\exp\left(\lambda_+x_3 - \lambda_- y_3 -\tau s -iy'\cdot\xi'\right)\,d\xi'd\,\tau\\
& = & (2\pi)^{-3} \int_{-\infty-i\sigma}^{\infty-i\sigma}\int_{\mathbb R^2}\exp\left(i\eta(t-s) + i(x'-y')\cdot\xi'\right) g(x,\,\xi',\,\eta;\,y_3,\,s)\,d\xi'\,d\eta,
\end{eqnarray*}
where
\begin{equation}\label{assume_g}
g(x,\,\xi',\,\eta;\,y_3,\,s)=\left[(-A_1 + A_2) \exp\left(\lambda_+x_3 - \lambda_- y_3\right)\right]_{\tau=i\eta}.
\end{equation}
Recall that Lemma \ref{Arima} still holds when the amplitude $g$ depends on $(x',\,y_3,\,s)$. So, if $g$ satisfies the assumptions of
Lemma \ref{Arima} with $\kappa=-1$, then we have the desired estimate for $x^\prime,\,y^\prime\in \mathbb R^2,\,x_3,\,y_3\leq 0,\,x_3-y_3\geq 0,\,t,\,s\in (0,\,T),\, t>s$ that
\begin{equation}\label{estimateOfLocalParametrix1}
|\dot G^+_{1(1)}(x,\,t;\,y,\,s)|\le c_1
(t-s)^{-\frac{3}{2}}\exp\Big(-c_2\frac{|x'-y'|^2+(x_3+y_3)^2}{t-s}\Big)
\end{equation}
with some positive constants $c_1$ and $c_2$, where we evaluated the inverse Laplace transform at $t-s$. In a similar way, we have
\begin{equation}\label{estimate_G_x1}
\left\{
\begin{array}{ll}
|\nabla_x \dot G^+_{1(1)}(x,\,t;\,y,\,s)|\le c_1
(t-s)^{-2}\exp\Big(-c_2\displaystyle\frac{|x'-y'|^2+(x_3+y_3)^2}{t-s}\Big),\\
\\
|\dot G^+_{1(2)}(x,\,t;\,y,\,s)|\le c_1
(t-s)^{-\frac{3}{2}}\exp\Big(-c_2\displaystyle\frac{|x'-y'|^2+(x_3-y_3)^2}{t-s}\Big),\\
\\
|\nabla_x \dot G^+_{1(2)}(x,\,t;\,y,\,s)|\le c_1
(t-s)^{-2}\exp\Big(-c_2\displaystyle\frac{|x'-y'|^2+(x_3-y_3)^2}{t-s}\Big)
\end{array} \right.
\end{equation}
for $x,\,y,\,t,\,s$ and $c_1,\,c_2$ as above.

\begin{remark}\label{estimate of principal part of parametrix}
From \eqref{estimateOfLocalParametrix1}, $\dot G^+_{1(1)}(x,\,t;\,y,\,s)\rightarrow 0$ as $t\downarrow s$ if $y_3<0$. We remark here again that the second and third estimates of \eqref{estimate_G_x1} are the estimates for the principal part of the fundamental solution of $\partial_t-\Delta$ in terms of the coordinates introduced by the coordinates transformation $\Phi_{x_0}$ flattening $\partial D$ near $x_0\in\partial D$.
\end{remark}

As an immediate consequence, we have the estimates
\begin{eqnarray}
&&|\dot G^+_{1,N}(x,\,t;\,y,\,s)| \leq c_1 (t-s)^{-\frac{3}{2}}\exp\Big(-c_2\frac{|x'-y'|^2+(x_3-y_3)^2}{t-s}\Big), \label{estimateOfLocalParametrix}\\
&&|\nabla_x \dot G^+_{1,N}(x,\,t;\,y,\,s)| \leq c_1 (t-s)^{-2}\exp\Big(-c_2\frac{|x'-y'|^2+(x_3-y_3)^2}{t-s}\Big) \label{estimate_G_x}
\end{eqnarray}
for $x^\prime,\,y^\prime\in \mathbb R^2,\,x_3,\,y_3\leq 0,\,x_3-y_3\geq 0,\,t,\,s\in (0,\,T),\, t>s$, where $c_1$ and $c_2$ are some positive constants. The same estimates hold for $G^-_{\ell,N}$ and $H^\pm_{\ell,N}$. If $t>s$, these are smooth with respect to $x,\,y,\,t,\,s$ by Remark \ref{Arima_R}.

\medskip
In the following, we first show that $\dot G_1^+(x,\,t;\,y,\,s)$ vanishes for $t<s$. Change the contour $\{ \tau=\sigma + i\mu:\, \mu\in\mathbb R \}$ of the integration with respect to $\tau$ for $\dot G_1^+(x,\,t;\,y,\,s)$ with $t<s$ to an infinitely large half circle $C_\infty:=\lim _{\rho\to\infty} C_\rho$ with
$C_{\rho}=\{ \sigma+\rho e^{i\theta}:\; -\frac{\pi}{2}\leq \theta \leq \frac{\pi}{2} \}.$ Then we can easily see that
\begin{equation}\label{p1}
\dot G_1^+(x,\,t;\,y,\,s)=0 \quad \mathrm{for}\;t<s.
\end{equation}

We next show that the assumptions of Lemma \ref{Arima} are satisfied for $g(\xi',\,\eta,\,\rho)$ with $\rho=|x_3-y_3|$. Let us verify the analyticity assumption. Then the other assumption can be easily justified. By \eqref{roots}, \eqref{coefficients} and \eqref{Ord-1Amplitude}, it is enough to prove that $\lambda_{\pm},\,\mu_{\pm}$ and their terms with square root are analytic, and for $y_3$ close to $0$
\begin{equation}\label{add1}
k(i\sum_{j=1}^2m_{3j}^0\xi_j+\mu_{+}m_{33}^0)-(i\sum_{j=1}^2m_{3j}^0\xi_j + \lambda_{+} m_{33}^0)
\end{equation}
does not vanish for $(\xi',\,\eta)\in L_\mu^2$. Here we note that $x_3<0$ is confined to $(-\delta, \, 0)$ with small $\delta>0$, so if $y_3\le-2\delta$ it is easy to see that $\dot G_1^+$ becomes a smoothing operator. This is why we can assume that $y_3$ is close to $0$.

Consider the characteristic equation for the operator $p_{2,0}^{(0)}$, that is,
\begin{equation*}
p_0(x')\xi_3^2 + p_1(x',\,\xi')\xi_3+(p_2(x',\,\xi')-\tau)=0,
\end{equation*}
where $$p_0(x')=-m_{33},\quad p_1(x',\,\xi')=-2\sum_{j=1}^2 m_{3j}\xi_j, \quad p_2(x',\,\xi')=-\sum_{i,j=1}^2 m_{ij}\xi_i\xi_j,\quad \tau=i\eta.$$
Its roots are given by
\begin{equation*}
\xi_3=\displaystyle\frac{-p_1 \pm z^{\pm}}{2p_0} \qquad \textrm{with }\,z^{\pm}=\sqrt{p_1^2-4p_0(p_2-\tau)},\;\pm \mathrm{Im}\, z^{\pm}>0.
\end{equation*}
By the ellipticity, there exists a constant $c>0$ such that $p_1^2-4p_0 p_2<-c|\xi'|^2$ for $\xi'\in\mathbb R^2\setminus\{0\},\,x'\in U$,
where $U$ is a bounded open set in which $p_0,\,p_1,\,p_2$ are smooth. Then we have the following claim:

\medskip
\begin{claim}\label{claim1}
There exists $\mu>0$ such that $p_1^2-4p_0 p_2 + 4p_0 i\eta\not\in [0,\,\infty)$ for $(\xi',\,\eta)\in L_{\mu}^{2}$.
\end{claim}

{\bf Proof.} We prove the claim by a contradiction argument. Note that
\begin{equation}\label{remind}
L_{\mu}^{2}\ni(\xi',\eta)\Longleftrightarrow \mathrm{Im}\,\eta < \mu(|\mathrm{Re}\, \eta|+|\mathrm{Re}\, \xi'|^2)-\mu^{-1}|\mathrm{Im}\, \xi'|^2.
\end{equation}
Suppose that for $(\xi',\,\eta)\in L_{\mu}^{2}$ there is a positive constant $m$ such that $p_1^2-4p_0 p_2 +4p_0 i\eta=m$. Set
$$\alpha:=-m_{33},\; \beta=(\beta_1,\,\beta_2):=-2(m_{31},\,m_{32}),\; \gamma=(\gamma_{ij})_{i,j=1,2}:=(-m_{ij})_{i,j=1,2}.$$
Then we obtain that
\begin{align*}
& p_0(x')=\alpha<0,\quad p_1(x',\,\xi')=\beta(x')\cdot\xi',\\
& p_2(x',\,\xi')=(\gamma(x')\xi')\cdot\xi'=(\gamma\xi')\cdot\xi'<0\quad \mbox{for}\; \xi'\in\mathbb R^{2}\setminus\{0\},\\
& p_1^2=\sum_{j,k=1}^2\beta_j \beta_k \xi_j \xi_k=(\beta\otimes\beta):(\xi'\otimes\xi'),\\
& m=(\beta\otimes\beta):(\xi^\prime\otimes\xi')-4\alpha(\gamma\xi')\cdot\xi'+4i\alpha\eta.
\end{align*}
For simplicity of notations, we denote $\xi'_R=\mathrm{Re}\,\xi',\,\xi'_I=\mathrm{Im}\,\xi',\,\eta_R=\mathrm{Re}\,\eta$ and $\eta_I=\mathrm{Im}\,\eta$. Note that
\begin{eqnarray*}
m & = &(\beta\otimes\beta):(\xi'_R\otimes\xi'_R)-(\beta\otimes\beta):(\xi'_I\otimes\xi'_I)
+2i(\beta\otimes\beta):(\xi'_R\otimes\xi'_I)\\
& & -4\alpha\{(\gamma\xi'_R)\cdot\xi'_R-(\gamma\xi'_I)\cdot\xi'_I\}-8i\alpha(\gamma\xi'_R)\cdot\xi'_I+4i\alpha\eta,
\end{eqnarray*}
which yields
\begin{eqnarray}
&&(\beta\otimes\beta):(\xi'_R\otimes\xi'_R)-(\beta\otimes\beta):(\xi'_I\otimes\xi'_I)
-4\alpha\{(\gamma\xi'_R)\cdot\xi'_R-(\gamma\xi'_I)\cdot\xi'_I\}-4\alpha\eta_I=m,\label{system_1}\\
&&(\beta\otimes\beta):(\xi'_R\otimes\xi'_I)-4\alpha(\gamma\xi'_R)\cdot\xi'_I+2\alpha\eta_R=0. \label{system_2}
\end{eqnarray}
We obtain from \eqref{system_1} that
\begin{equation}\label{cl}
\{-(\beta\otimes\beta):(\xi'_R\otimes\xi'_R)+4\alpha(\gamma\xi'_R)\cdot\xi'_R\}+\{(\beta\otimes\beta):(\xi'_I\otimes\xi'_I)
-4\alpha(\gamma\xi'_I)\cdot\xi'_I\}=-m-4\alpha\eta_I.
\end{equation}
The left hand side (LHS) of \eqref{cl} has the estimate that $\mathrm{LHS}> c'|\xi'_R|^2-c''|\xi'_I|^2$ for some positive constants $c'$ and
$c''$. For the right hand side (RHS) of \eqref{cl}, by the definition of $L_{\mu}^2$, we deduce from \eqref{system_2} that
\begin{eqnarray*}
\mathrm{RHS} &\leq& -4\alpha\eta_I
<(-4\alpha)\{\mu(|\eta_R|+|\xi'_R|^2)-\mu^{-1}|\xi'_I|^2\}\\
&=&(-4\alpha)\{\mu((-2\alpha)^{-1}|(\beta\otimes\beta):(\xi'_R\otimes\xi'_R)-4\alpha(\gamma\xi'_R)\cdot\xi'_I|+|\xi'_R|^2)-\mu^{-1}|\xi'_I|^2\}\\
&\leq& \mu K(|\xi'_R|^2+|\xi'_I|^2)-(-4\alpha)\mu^{-1}|\xi'_I|^2
\end{eqnarray*}
for some positive constant $K$. Thus, we have
\begin{equation*}
\tilde{c}(|\xi'_R|^2+|\xi'_I|^2)\leq \mu K (|\xi'_R|^2+|\xi'_I|^2), \qquad \mu>0.
\end{equation*}
By taking $\mu$ sufficiently small, we have $\xi'_R=\xi'_I=0$ and hence $\eta_R=0$ by \eqref{system_2}. Then, $-4\alpha\eta_I=m$ gives
$\eta_I\geq 0$. This contradicts to \eqref{remind}. The proof of this claim is now complete.  \hfill $\Box$

\medskip
Using the above claim, we can easily see that $\lambda_{\pm}$, $\mu_{\pm}$ and their terms with square root are analytic for $(\xi',\,\eta)\in L_\mu^2$. In addition, for  $y_3$ near to $0$, there exists a constant $\mu>0$ such that \eqref{add1} does not vanish. We give the proof as follows. Define
\begin{equation*}
R_0:=\sum_{j=1}^2m_{3j}^0\xi_j,\qquad R_1:=\sum_{j=1}^2m_{3j}^1\xi_j,\qquad Q:=\sum_{j,l=1}^2m_{jl}^1\xi_j\xi_l.
\end{equation*}
Recalling the definitions of $\lambda_{\pm}$ and $\mu_{\pm}$ in \eqref{roots}, we have
\begin{eqnarray*}
& &k(i\sum_{j=1}^2m_{3j}^0\xi_j+\mu_{+}m_{33}^0)-(i\sum_{j=1}^2m_{3j}^0\xi_j + \lambda_{+} m_{33}^0)\\
&=&k(iR_0 + \mu_+ m^0_{33}) - (iR_0 + \lambda_+ m^0_{33})\\
&=&i(k-1)\left(R_0 - m^0_{33} (m^1_{33})^{-1} R_1 \right) \\
& & +  m^0_{33} (m^1_{33})^{-1} \left(\sqrt{m^1_{33}(k^2Q+k\tau) - k^2 R_1^2} - \sqrt{m^1_{33}(Q+\tau) - R_1^2} \right).
\end{eqnarray*}
When $y_3$ is near to $0$, $R_0 - m^0_{33} (m^1_{33})^{-1} R_1$ is also near to $0$. By the same argument as that for Claim \ref{claim1}, there exists a constant $\mu>0$ such that
\begin{equation*}
(k-1)[m^1_{33}(k+1)Q+m^1_{33}\tau -(k+1)R_1^2] \not\in [0,\,\infty) \qquad  \mathrm{for}\; (\xi',\,\eta)\in L_{\mu}^{2},
\end{equation*}
which implies that
$$\sqrt{m^1_{33}(k^2Q+k\tau) - k^2 R_1^2} - \sqrt{m^1_{33}(Q+\tau) - R_1^2}\not = 0.$$
This completes the proof of showing the analyticity assumption for \eqref{assume_g}, and further by a scaling argument we can see that $g$ given by \eqref{assume_g} satisfies the estimate \eqref{eq:estimate for g} with $\kappa=-1$.

\bigskip
In conclusion, we have justified that $g$ defined by \eqref{assume_g} satisfies the assumptions of Lemma \ref{Arima}, and hence we obtain the desired estimate \eqref{estimateOfLocalParametrix} for $\dot G_1^+$ where we only take the leading term $\mathbf a_{1,-1}$ of the amplitude $\mathbf a_1$. In the following, we estimate the error terms if we truncate the amplitudes at $j=N$. That is to consider $\sum_{j=0}^N \mathbf a_{1,1-j}$ for $\mathbf a_1$. Let
$$\dot G_{1,N}(x,\,t;\,y,\,s):=\dot G_{1,N}^\pm(x,\,t;\,y,\,s), \qquad \pm(x_3-y_3)>0$$
(see the paragraph just after \eqref{w1} for the definition of $\dot G_{1,N}^\pm(x,\,t;\,y,\,s)$). Then, by the construction of amplitudes and Lemma \ref{Arima}, the error
\begin{equation*}
R_{1,N}:=R_{1,N}(x,\,t;\,y,\,s)=\left[\partial_t \dot G_{1,N} - \mathcal J^{-1}(x)\nabla_x \cdot \big(\mathcal J(x)M(x)\nabla_x \dot G_{1,N}\big)\right] - I
\end{equation*}
is smooth enough for $x_3\leq 0$ with respect to all the variables. Moreover, by the expression \eqref{A_L} of $\mathcal A_{\ell,1-L}$, we see that $\mathcal A_{1,1}=\mathcal A_{1,0}=\cdots=\mathcal A_{1,1-N}=0$ and $\mathcal A_{1,-N}$ is the dominant term in the remaining part of $\mathcal A_1$ with $\mathrm{ord}\,\mathcal A_{1,-N}=-N$. Hence, using the estimate \eqref{eq:estimate for g} with $\kappa=-N$, we have the estimate
\begin{equation}\label{error estimate}
|\partial_t^j \partial_x^\alpha R_{1,N}(x,\,t;\,y,\,s)|\leq
C_{N,M}'(t-s)^{-2+\frac{N-|\alpha|}{2}-j}\exp\Big(-C_{N,M}\frac{|x-y|^2}{t-s}\Big)
\end{equation}
for any $j\in{\mathbb Z}_+,\,\,\alpha\in{\mathbb Z}_+^3$ such that $2j+|\alpha|\le M$ with $M\in \mathbb N$ and $M\le N-4$, where $C_{N,M}$ and $C_{N,M}'$ are positive constants.

We define $\dot H_{1,N}(x,\,t;\,y,\,s),\,\dot G_{2,N}(x,\,t;\,y,\,s),\,\dot H_{2,N}(x,\,t;\,y,\,s)$ and the associated error terms $S_{1,N}(x,\,t;\,y,\,s)$, $R_{2,N}(x,\,t;\,y,\,s)$, $S_{2,N}(x,\,t;\,y,\,s)$ analogously to $\dot G_{1,N}(x,\,t;\,y,\,s)$ and $R_{1,N}(x,\,t;\,y,\,s)$, respectively. They satisfy similar properties as those of $\dot G_{1,N}$ and $R_{1,N}$. Since the pairs $(\dot G_{\ell,N}, \dot H_{\ell,N})$ for $\ell=1,\,2$ satisfy the boundary condition of \eqref{eq:L_ITP-2} locally in terms of the coordinates $\xi=(\xi_1,\,\xi_2,\,\xi_3)=\Phi_{x_0}(x)$, these pairs can be used to define a local parametrix for \eqref{eq:ITP-2} in the open neighborhood $U(x_0)$ of $x_0\in\partial D$.

\section{Construction of the Green function for parabolic ITP}\label{patch}
\setcounter{equation}{0}

In this section, using a partition of unity, we patch the local parametrices constructed above and the fundamental solution of the diffusion equation so that we have a global parametrix for \eqref{eq:ITP-2}. Then, using the Levi method, we construct the Green function from this parametrix.

Take $x_0^{(j)}\in \partial D\,(j=1,\,2,\,\cdots,\,J)$ so that $\{U_j:=U(x_0^{(j)})\}_{j=1}^J$ is an open covering of $\partial D$. Let $U_0$ be an open set such that $\overline{U_0}\subset D$ and $\{U_j\}_{j=0}^J$ gives an open covering of $\overline D$. Let $\varphi_j\in C_0^\infty (U_j)\,(j=0,\,1,\,\cdots,\,J)$ be a partition of unity subordinated to this cover and $\psi_j\in C_0^\infty (U_j)\,(j=0,\,1,\,\cdots,\,J)$ satisfy $\psi_j=1$ on $\mathrm{supp}\,\varphi_j\,(j=0,\,1,\,\cdots,\,J)$.

We will abuse the notation $(G_\ell^j,\,H_\ell^j)$ to denote the local parametrix $\big(\dot G_{\ell,N}(x,t;y,s),\,\dot H_{\ell,N}(x,t;y,s)\big)$ constructed for each $U_j:=U(x_0^{(j)})$ $(j=1,2,\cdots,J)$ in the previous section.

Let us first look for a parametrix for $(G_1,\,H_1)$. Set
\begin{equation}\label{fundamental solutions}
\left\{
\begin{array}{lll}
G_1^0(x,\,t;\,y,\,s)&:=&H(t-s)\,(4\pi(t-s))^{-3/2}\,\exp\left(-(4(t-s))^{-1} |x-y|^2\right),\\
H_1^0(x,\,t;\,y,\,s)&:=&H(t-s)\,(4k\pi(t-s))^{-3/2}\,\exp\left(-(4k(t-s))^{-1} |x-y|^2\right),
\end{array}
\right.
\end{equation}
where $H$ is the Heaviside function defined by $$H(t-s)=\begin{cases} 1, & t-s>0,\\ 0, &t-s\leq 0. \end{cases}$$
Since the Green function ${\mathbb G}$ (see \eqref{Green function G}) is a distribution belonging to ${\mathscr D}^\prime(D_T\times D_T)$ and its singularities are only near the diagonal $\{(x,\,t;\,y,\, s)\in D_T\times D_T:\, x=y,\, t=s\}$, we define a parametrix $ (\tilde G_1^\prime(x,\,t;\,y,\,s),\,\tilde H_1^\prime(x,\,t;\,y,\,s))$ for $K_1(x,\,t;\,y,\,s)$ (see just after \eqref{F ell}) by
\begin{equation}\label{parametrices tilde G_1 prime}
\left\{\begin{array}{lll}
\tilde G_1^\prime(x,\,t;\,y,\,s)&:=&\displaystyle\sum_{j=0}^J\psi_j(x)G_1^j(x,\,t;\,y,\,s)\varphi_j(y),\\
\tilde H_1^\prime(x,\,t;\,y,\,s)&:=&\displaystyle\sum_{j=0}^J\psi_j(x)H_1^j(x,\,t;\,y,\,s)\varphi_j(y).
\end{array} \right.
\end{equation}
For fixed $s\in[0,\,T]$, let $\tilde G_1^\prime (t,\,s),\,\tilde H_1^\prime(t,\,s),\,G_1^j(t,\,s),\,H_1^j(t,\,s)$ for $j=0,\,1,\,\cdots,\,J$ be the pseudo-differential operators with parameter $t\in[0,\,T]$ whose Schwartz kernels are
$\tilde G_1^\prime (x,\,t;\,y,\,s),\, \tilde H_1^\prime(x,\,t;\,y,\,s),\, G_1^j(x,\,t;\,y,\,s),$ $H_1^j(x,\,t;\,y,\,s)$, respectively. For example, $(\tilde G_1^\prime(t,\,s)\phi)(x)$ is defined as
\begin{equation}\label{tilde G1 prime}
\left(\tilde G_1^\prime(t,\,s)\phi\right)(x):=\int_D \tilde G_1^\prime(x,\,t;\,y,\,s)\phi(y)\,dy
\end{equation}
for any function $\phi(x)\in C_0^\infty(D)$. Sometimes we suppress the parameter $t$ and denote $\tilde G_1^\prime (t,\,s),\,\tilde H_1^\prime(t,\,s),$ $G_1^j(t,\,s),\,H_1^j(t,\,s)$ by $\tilde G_1^\prime (s),\,\tilde H_1^\prime(s),\,G_1^j(s),\,H_1^j(s)$, respectively.

Note that
\begin{eqnarray*}
(\partial_t - \Delta)\tilde G_1^\prime(s) &=& \sum_{j=0}^J \left\{ \psi_j (\partial_t-\Delta)G_1^j(s)\varphi_j - [\Delta,\,\psi_j] G_1^j(s)\varphi_j \right\}\\
&=:&I + S_{\tilde G_1^\prime}^N(s) - \sum_{j=0}^J  [\Delta,\,\psi_j] G_1^j(s)\varphi_j\\
&=:&I + \dot R _{\tilde G_1^\prime}^N(s),
\end{eqnarray*}
where $[\Delta,\,\psi_j]:=\nabla\psi_j\cdot\nabla - \Delta\psi_j$ is the commutator of $\Delta$ and the multiplication by $\psi_j$. Using the estimate \eqref{error estimate} and its derivation, we know that
$$S_{\tilde G_1^\prime}^N(s):\,C^m([0,\,T];\,H^r(D))\to C^m([0,\,T];\,H^{r+m}(D))\qquad \textrm{with }\, m\in \mathbb Z_+,\; 3m\le N-7$$
is a bounded operator for any $r\in{\mathbb Z}_+$ vanishing for $t<s$ and vanishes at $t=s$ by order $m$. In the sequel, we fix $m$ as above.

Since $\mathrm{supp}\,\varphi_j \cap \mathrm{supp}\,[\Delta,\,\psi_j]=\emptyset$, it can be seen that $\sum_{j=0}^J[\Delta,\,\psi_j]G_1^j(s)\varphi_j$ is a smoothing operator flat at $t=s$. Hence, $\dot R_{\tilde G_1^\prime}^N(s):\,C^m([0,\,T];H^r(D))\to C^m([0,\,T];\,H^{r+m}(D))$ is a bounded operator for any $r\in{\mathbb Z}_+$ vanishing for $t<s$ and vanishes at $t=s$ by order $m$. Similarly, we can derive that
\begin{equation*}
(\partial_t - k\Delta)\tilde H_1^\prime(s)=\dot R _{\tilde H_1^\prime}^N(s),
\end{equation*}
where $\dot R_{\tilde H_1^\prime}^N(s):\,C^m([0,\,T];\,H^r(D))\to C^m([0,\,T];\,H^{r+m}(D))$ is a bounded operator for any $r\in{\mathbb Z}_+$ vanishing for $t<s$ and vanishes at $t=s$ by order $m$ which can be shown likewise \eqref{p1}. From the construction of the local parametrix $(G_1^j,\,H_1^j)$, it is clear that
\begin{equation*}
\tilde H_1^\prime(s) - \tilde G_1^\prime(s) = 0 \quad\mathrm{on}\;\,\partial D.
\end{equation*}
Further, since $\mathrm{supp}\,\varphi_j \cap \mathrm{supp}\,[\partial_\nu,\,\psi_j]=\emptyset$, we have
\begin{eqnarray*}
k \partial_\nu \tilde H_1^\prime(s) - \partial_\nu \tilde G_1^\prime(s) &=& \sum_{j=1}^J \left\{ \psi_j (k\partial_\nu H_1^j(s) - \partial_\nu G_1^j(s))\varphi_j - (k[\partial_\nu,\,\psi_j]H_1^j(s) - [\partial_\nu,\,\psi_j]G_1^j(s))\varphi_j  \right\}\\
&=:& \Lambda(s)
\end{eqnarray*}
with a bounded smoothing operator $\Lambda(t,\,s)=\Lambda(s)$ flat at $t=s$, where we have used the previous convention given just after \eqref{tilde G1 prime}. By noticing $G_1^j(t,\,s)|_{t<s}=H_1^j(t,\,s)|_{t<s}=0$ for $j=0,\,1,\,\cdots,\,J$ from their definitions, we have
\begin{equation*}
\tilde G_1^\prime(t,\,s)|_{t<s}=\tilde H_1^\prime(t,\,s)|_{t<s}=0.
\end{equation*}
In the sequel, we will use for example the notation $\tilde G_1^\prime(t,\,s)$ for $\tilde G_1^\prime(s)$ without mentioning the convention if it is clear from the context.

Finally, by moving $\Lambda(t,\,s)$ to inhomogeneous terms via the inverse trace operator, there exists a parametrix $(\tilde G_1(t,\,s),\,\tilde H_1(t,\,s))$ which has the same estimates as those of $(\tilde G_1^\prime(t,\,s),\,\tilde H_1^\prime(t,\,s))$ and satisfies
\begin{equation}\label{eq:ITP-p}
\begin{cases}
(\partial_t-\Delta)\tilde G_1(t,\,s)=I+R_{G_1}^N(t,\,s) & \textrm{ in }D_T,\\
(\partial_t-k\Delta)\tilde H_1(t,\,s)=R_{H_1}^N(t,\,s) & \textrm{ in }D_T,\\
\tilde G_1(t,\,s) - \tilde H_1(t,\,s)=0 & \textrm{ on }(\partial D)_T,\\
\partial_{\nu}\tilde G_1(t,\,s) - k\partial_{\nu}\tilde H_1(t,\,s) =0 & \textrm{ on }(\partial D)_T, \\
\tilde G_1(t,\,s)=\tilde H_1(t,\,s)=0 & \textrm{ for } t<s,
\end{cases}
\end{equation}
where $R_{G_1}^N(t,\,s),\,R_{H_1}^N(t,\,s):\,C^m([0,\,T];\,H^r(D))\to C^m([0,\,T];\,H^{r+m}(D))$ are bounded operators for any $r\in{\mathbb Z}_+$ vanishing for $t<s$ and vanish at $t=s$ by order $m$.

In the same way, we can show that there exists a parametrix $(\tilde G_2(t,\,s),\,\tilde H_2(t,\,s))$ such that
\begin{equation}\label{eq:ITP-p2}
\begin{cases}
(\partial_t-\Delta)\tilde G_2(t,\,s)=R_{G_2}^N(t,\,s) & \textrm{ in }D_T,\\
(\partial_t-k\Delta)\tilde H_2(t,\,s)=I+R_{H_2}^N(t,\,s) & \textrm{ in }D_T,\\
\tilde G_2(t,\,s)-\tilde H_2(t,\,s)=0 & \textrm{ on }(\partial D)_T,\\
\partial_{\nu}\tilde G_2(t,\,s)-k\partial_{\nu}\tilde H_2(t,\,s)=0 & \textrm{ on }(\partial D)_T, \\
\tilde G_2(t,\,s)=\tilde H_2(t,\,s)=0 & \textrm{ for } t<s,
\end{cases}
\end{equation}
where $R_{G_2}^N(t,\,s),\,R_{H_2}^N(t,\,s):\,C^m([0,\,T];\,H^r(D))\to C^m([0,\,T];\,H^{r+m}(D))$ are bounded operators for any $r\in{\mathbb Z}_+$ vanishing for $t<s$ and vanish at $t=s$ by order $m$. Thus, we have a matrix
\begin{equation*}
\tilde {\mathbb G}(t,\,s)=\left(
\begin{array}{lr}
\tilde G_1(t,\,s) & \tilde G_2(t,\,s) \\
\tilde H_1(t,\,s) & \tilde H_2(t,\,s)
\end{array} \right)
\end{equation*}
such that
\begin{equation*}
\mathcal L \tilde {\mathbb G} = I + R_N, \qquad  \tilde {\mathbb G}(t,\,s)=0\;\, \mathrm{for}\;\, t<s,
\end{equation*}
where
\begin{equation*}
\mathcal L = \left(
\begin{array}{lr}
\partial_t - \Delta & 0 \\
0 & \partial_t - k \Delta
\end{array} \right),
\qquad
R_N(t,\,s) = \left(
\begin{array}{lr}
R_{G_1}^N(t,\,s) & R_{G_2}^N(t,\,s) \\
R_{H_1}^N(t,\,s) & R_{H_2}^N(t,\,s)
\end{array} \right)
\end{equation*}
and $R_N(t,\,s): C^m([0,\,T];\,H^r(D))\to C^m([0,\,T];\,H^{r+m}(D))$ is a bounded operator for any $r\in{\mathbb Z}_+$ vanishing for $t<s$ and vanishes at $t=s$ by order $m$. Further, $\tilde {\mathbb G}(t,\,s)$ satisfies the following estimates. That is, for any $x,\,y\in D,\,t,\,s\in(0,\,T),\,t>s$, the estimates
\begin{eqnarray}\label{estimate of tilde mathbb G}
&&\tilde{|\mathbb G}(x,\,t;\,y,\,s)| \leq c_1 (t-s)^{-3/2}\exp\Big(-c_2\frac{|x-y|^2}{t-s}\Big)\label{estimate of tilde mathbb G 1}, \\
&&|\nabla_x \tilde{\mathbb G}(x,\,t;\,y,\,s)| \leq c_1 (t-s)^{-2}\exp\Big(-c_2\frac{|x-y|^2}{t-s}\Big)
\end{eqnarray}
hold with positive constants $c_1$ and $c_2$ independent of $x,\,y,\,t,\,s$.

Define $\tilde {\mathcal G}(t,\,s)$ by
\begin{equation}\label{parametrix for Cauchy prob}
\tilde {\mathcal G}(t,\,s)=\tilde {\mathbb G}(t,\,s)\big|_{t\ge s}.
\end{equation}
Namely, $\tilde {\mathcal G}(t,\,s)$ is the restriction of $\tilde {\mathbb G}(t,\,s)$ to $t \geq s$. Then $\tilde {\mathcal G}$ is a parametrix for the following Cauchy problem:
\begin{equation}\label{Cauchy_0}
\begin{cases}
(\partial_t-\Delta)v=0 & \textrm{ in }D\times (s,\,T),\\
(\partial_t-k\Delta)u=0 & \textrm{ in }D\times (s,\,T),\\
v-u=0 & \textrm{ on }(\partial D)\times (s,\,T),\\
\partial_{\nu}v - k\partial_{\nu}u=0 & \textrm{ on }(\partial D)\times (s,\,T), \\
v=v_0,\; u=u_0& \textrm{ at } t=s
\end{cases}
\end{equation}
with
\begin{equation}\label{parametrix for Cauchy problem}
\mathcal L \tilde {\mathcal G}=R_N,\quad \tilde {\mathcal G}(s,\,s)=I.
\end{equation}

In order to see the second equation of \eqref{parametrix for Cauchy problem},  denote $\tilde{\mathcal{G}}^+(t,\,s)=\tilde{\mathcal{G}}(t\,,s)$ for $t\ge s$ and $\tilde{\mathcal{G}}^-(t,\,s)=0$ for $t<s$. Then $\tilde{\mathbb{G}}(t,\,s)$ is of course given as
\begin{equation*}
\tilde{\mathbb G}(t,\,s):=
\begin{cases}
\tilde{\mathcal G}^+(t,\,s), & t-s\ge 0,\\
\tilde{\mathcal G}^-(t,\,s), & t-s<0.
\end{cases}
\end{equation*}
Using the expression of the parametrix $\tilde{\mathcal G}^+(t,\,s)=\tilde{\mathcal{G}}(t,\,s)$, its equation and \eqref{estimate of tilde mathbb G 1}, $\tilde{\mathcal G}^+ \in C^1\big(\{\pm(t-s)\ge0\}\cap (0,\,T),\,\mathscr D^\prime(D)\big)$.
For any $\varphi\in C_0^\infty(D_T)$ and any fixed $(y,\,s)\in D_T$, we have
\begin{eqnarray*}
&&\langle\delta(\cdot-y)\delta(\cdot-s)I,\,\varphi\rangle + \langle R_N,\varphi\rangle\\
&=&\langle (\partial_t - L)\tilde{\mathbb G},\,\varphi\rangle\\ &=& -\langle \tilde{\mathbb G},\,\partial_t \varphi\rangle - \langle L \tilde{\mathbb G},\,\varphi\rangle\\
&=& -\langle\tilde{\mathcal G}^+,\,\partial_t \varphi\rangle -\langle\tilde{\mathcal G}^-,\,\partial_t \varphi\rangle - \langle L \tilde{\mathbb G},\,\varphi\rangle\\
&=& -\lim_{\varepsilon \to 0} \int_{s+\varepsilon}^T \langle \tilde{\mathcal G}^+,\,\partial_t \varphi(t,\,\cdot) \rangle\, dt
-\lim_{\varepsilon \to 0} \int_0^{s-\varepsilon} \langle \tilde{\mathcal G}^-,\,\partial_t \varphi(t,\,\cdot) \rangle\, dt - \langle L \tilde{\mathcal G}^+,\,\varphi\rangle\\
&=& \langle \tilde{\mathcal G}^+|_{t=s},\, \varphi(s,\,\cdot)\rangle + \langle \partial_t \tilde{\mathcal G}^+,\,\varphi\rangle - \langle L \tilde{\mathcal G}^+,\,\varphi\rangle\\
&=& \langle \tilde{\mathcal G}^+|_{t=s},\, \varphi(s,\,\cdot)\rangle +\langle R_N,\varphi\rangle,
\end{eqnarray*}
where $L=\partial_t-\mathcal L$. Thus we have $\tilde {\mathcal G}(s,\,s)=I$. An alternative direct proof of this will be given in Appendix 2.

In the following, we construct the Green function $\mathcal G$ for \eqref{Cauchy_0} from the parametrix $\tilde {\mathcal G}$ by using the Levi method. Fixing $r\in{\mathbb Z}_+$, let
\begin{eqnarray*}
W_1(t,\,s)&=&-R_N(t,\,s),\\
W_j(t,\,s)&=&\int_s^t W_1(t,\,s')\, W_{j-1}(s',\,s)\,ds', \quad j\geq 2.
\end{eqnarray*}
Here we note that for instance the operator $R_N(t,\,s)$ should be understood as an
integral operator on $C^m([s,\,T];\,H^r(D))$ with kernel $R_N(x,\,t;\,y,\,s)$, that is,
$$
(R_N(t,\,s)\phi)(x):=\int_s^t\int_D  R_N(x,\,t;\,y,\,s')\phi(y,\,s')\,dy\,ds'.
$$
We have
\begin{equation}\label{add2}
\sum_{j=1}^l W_j(t,\,s) = -R_N(t,\,s) - \int_s^t R_N(t,\,s') \sum_{j=1}^{l-1}W_j(s',\,s)\,ds'.
\end{equation}
Notice that $R_N(t,\,s):\,C^m([s,\,T];\,H^r(D))\to C^m([s,\,T];\,H^{r}(D))$ is uniformly bounded for each $r\in{\mathbb Z}_+$ vanishing for $t<s$ and vanishes at $t=s$ by order $m$. Let $\Vert\cdot\Vert$ denote the operator norm for operators on $C^m([s,\,T];\,H^r(D))$. Then we have the following estimates:
\begin{eqnarray*}
\| W_1(t,\,s) \| &\leq & C_0,\\
\| W_j(t,\,s) \| &\leq & \frac{C_0^{j-1}}{(j-1)!}(t-s)^{j-1},\quad j\geq 2.
\end{eqnarray*}
Therefore, it can be easily seen that
\begin{equation*}
W(t,\,s):= \sum_{j=1}^\infty W_j(t,\,s)
\end{equation*}
converges as a bounded operator on $C^m([s,\,T];\,H^r(D))$ and vanishes for $t<s$. Furthermore, we observe from \eqref{add2} that
\begin{equation}\label{add3}
W(t,\,s) = -R_N(t,\,s) - \int_s^t R_N(t,\,s')W(s',\,s)\,ds'.
\end{equation}
Set
\begin{equation}\label{add4}
{\mathcal G} := \tilde{\mathcal G}+ \int_s^t \tilde{\mathcal G}(t,\,s')W(s',\,s)\,ds'.
\end{equation}
Note that by \eqref{estimateOfLocalParametrix}, \eqref{estimate_G_x}, \eqref{fundamental solutions} and the fact in Appendix 3, $\tilde{\mathcal G}(t,\,s):\,L^2(D_T)\to L^2((s,\,T);\,H^1(D))$ is bounded.
By the first equation of \eqref{parametrix for Cauchy problem}, \eqref{add3}, the direct calculations give that
\begin{eqnarray*}
\mathcal L {\mathcal G} &=& \mathcal L \tilde {\mathcal G} + \mathcal L \left(\int_s^t  \tilde {\mathcal G}(t,\,s')W(s',\,s)\,ds' \right)\\
&=& R_N(t,\,s) + W(t,\,s) + \int_s^t R_N(t,\,s')W(s',\,s)\,ds'\\
&=& 0.
\end{eqnarray*}
Also from the second equation of \eqref{parametrix for Cauchy problem} and \eqref{add4}, we have
\begin{equation}\label{initial condition for tilde mathcal G}
\mathcal{G}(s,\,s)=I.
\end{equation}
Thus, we have completed our argument of constructing the Green function $\mathbb G(x,\,t;\,y,\,s)$ for the parabolic interior transmission problem \eqref{eq:ITP-2}, which is given by
\begin{equation}\label{G}
\mathbb G(x,\,t;\,y,\,s)=
\begin{cases}
\mathcal G (x,\,t;\,y,\,s),& t\geq s,\\
0, & t <s.
\end{cases}
\end{equation}

From the construction of the Green function $\mathbb G(x,\,t;\,y,\,s)$, we have for $x,\,y\in D,\,t,\,s\in(0,\,T),\,t>s$ that
\begin{eqnarray}\label{estimate of mathbb G}
&&|\mathbb G(x,\,t;\,y,\,s)| \leq c_1 (t-s)^{-3/2}\exp\Big(-c_2\frac{|x-y|^2}{t-s}\Big), \\
&&|\nabla_x \mathbb G(x,\,t;\,y,\,s)| \leq c_1 (t-s)^{-2}\exp\Big(-c_2\frac{|x-y|^2}{t-s}\Big),
\end{eqnarray}
where $c_1$ and $c_2$ are positive constants independent of $x,\,y,\,t,\,s$. Therefore, for any given inhomogeneous term $N=(N_1,\,N_2)\in L^2(D_T)$,  $U=(v,\,u)$ given by
\begin{equation*}
U(x,\,t) = \int_0^t \int_D \mathbb G(x,\,t;\,y,\,s)\, N(y,\,s)\,dyds.
\end{equation*}
is a solution of \eqref{eq:ITP-2}. Further, by the fact given in Appendix 3, this $U=U(x,\,t)$ belongs to $L^2((0,\,T);\, H^1(D))$.
In the next section, we will show the unique solvability of \eqref{eq:ITP-2}. The meaning of the unique solvability is as follows. Define the space of solutions $W((0,\,T))\subset L^2((0,\,T);\, H^1(D))$ to  \eqref{eq:ITP-2} by
\begin{equation}\label{solution space}
W((0,\,T))=\{ U(x,\,t)= (\mathbb G N)(x,\,t)=\int_0^T\int_D \mathbb G(x,\,t;\,y,\,s)\, N(y,\,s)\,dyds:\; N\in L^2(D_T)\}.
\end{equation}
Note that we can also write $U(x,\,t)$ in the form
\begin{equation}
U(x,\,t)=\int_0^t\int_D \mathcal{G}(x,\,t;\,y,\,s)\, N(y,\,s)\,dyds,
\end{equation}
which we abbreviate in the form
\begin{equation}
U(t)=\int_0^t \mathcal{G}(t,\,s)\, N(s)\,ds.
\end{equation}
Then we say that for any given $N\in L^2(D_T)$, (1.5) is uniquely solvable if there exist a unique solution in $W((0,\,T))$.

\section{Uniqueness of the Green function and solvability of ITP}\label{unique}
\setcounter{equation}{0}

In this section, we show the uniqueness of the Green function for the parabolic ITP \eqref{eq:ITP-2}, and then prove the unique solvability of \eqref{eq:ITP-2}. To this end, we use the duality argument given in \cite[Chapter 2]{Ito} and introduce the adjoint problem of the problem \eqref{Cauchy_0}:
\begin{equation}\label{Cauchy_2}
\begin{cases}
(-\partial_s-\Delta)w=0 & \textrm{ in }D\times (0,\,t),\\
(-\partial_s-k\Delta)z=0 & \textrm{ in }D\times (0,\,t),\\
w+z=0 & \textrm{ on }(\partial D)\times (0,\,t),\\
\partial_{\nu}w+k\partial_{\nu}z=0 & \textrm{ on }(\partial D)\times (0,\,t), \\
w=w_0,\; z=z_0& \textrm{ at } s=t.
\end{cases}
\end{equation}
We remark that the adjoint problem \eqref{Cauchy_2} also satisfies the Lopatinskii condition and the Green function can be constructed in the same way as for our original ITP. Denote by $\mathcal H(y,\,s;\,x,\,t)$ a Green function for \eqref{Cauchy_2}. Then it is enough to show that for any Green function $\mathcal G(x,\,t;\,y,\,s)$, we have
\begin{equation}\label{uq2}
\mathcal G(x,\,t;\,y,\,s) = \mathcal H^*(y,\,s;\,x,\,t),
\end{equation}
where $\mathcal H^*(y,\,s;\,x,\,t)$ is the adjoint of $\mathcal H(y,\,s;\,x,\,t)$.

For $\varpi<\vartheta$, define the function space $F((\varpi,\,\vartheta))$ by
$$F((\varpi,\,\vartheta)):=\left\{V:\,V\in H^1((\varpi,\,\vartheta);\,L^2(D))\cap L^2((\varpi,\,\vartheta);\,H^2(D))\right\}.$$
For $0<s<\tau<t<T$, if $U=(v,\,u)^{\mathbb T}\in F((s,\,T)),\,Z=(w,\,z)^{\mathbb T}\in F((0,\,t))$, where $(v,\,u)$ and $(w,\,z)$ are solutions to \eqref{Cauchy_0} and \eqref{Cauchy_2}, respectively, then we have
\begin{eqnarray}\label{uq1}
\int_D \partial_\tau(U\cdot Z)\,dx &=& \int_D\left[ (w\partial_\tau v + v\partial_\tau w) + (z\partial_\tau u + u \partial_\tau z) \right]dx \nonumber\\
&=& \int_D\left[ (w\Delta v - v\Delta w) + k(z\Delta u - u \Delta z) \right]dx \nonumber\\
&=& \int_{\partial D}\left[ (w\partial_\nu v - v\partial_\nu w) + k(z\partial_\nu u - u \partial_\nu z) \right]dx \nonumber\\
&=& \int_{\partial D}\left[ -k(z\partial_\nu u - u\partial_\nu z) + k(z\partial_\nu u - u \partial_\nu z) \right]dx \nonumber\\
&=& 0.
\end{eqnarray}
This means that
\begin{equation*}
\int_D U(z,\,\tau)\cdot Z(z,\,\tau)\,dz
\end{equation*}
is independent of $\tau$, and leads to
\begin{equation}\label{uq3}
\int_D U(z,\,t)\cdot Z(z,\,t)\,dz = \int_D U(z,\,s)\cdot Z(z,\,s)\,dz.
\end{equation}

\medskip
We now take $U_0,\,Z_0\in C_0^\infty (D)$ and set
\begin{eqnarray}
U(z,\,\tau)&=&\int_D \mathcal G(z,\,\tau;\,y,\,s)\,U_0(y)\,dy, \label{U}\\
Z(z,\,\tau)&=&\int_D \mathcal H(z,\,\tau;\,x,\,t)\,Z_0(x)\,dx. \label{Z}
\end{eqnarray}
Then we can prove that $U\in F((s,\,T))$ and $Z\in F((0,\,t))$. Here we only give the proof for $U\in F((s,\,T))$, since $Z\in F((0,\,t))$ can be proven in the same way.

For any fixed $s$, we simply write \eqref{U} as $U(\tau)=\mathcal G(\tau,\,s) U_0$. Then $\partial_\tau U(\tau)=\partial_\tau\mathcal G(\tau,\,s) U_0$. For any $U_0 \in C_0^\infty(D)$, we want to show $\partial_\tau \mathcal G(\tau,\,s)U_0 \in L^2((s,\,T);\,L^2(D))$ and $\partial_x^\alpha \mathcal G(\tau,\,s)U_0 \in L^2((s,\,T);\,L^2(D))$ for $|\alpha|\le2$. By \eqref{add4}, we have
\begin{eqnarray}\label{N1}
\partial_\tau \mathcal G(\tau,\,s) &=& \partial_\tau \tilde{\mathcal G}(\tau,\,s) + \tilde {\mathcal G}(\tau,\,\tau)W(\tau,\,s) + \int_s^\tau \partial_\tau \tilde{\mathcal G}(\tau,\,s^\prime)\, W(s^\prime,\,s)\, ds^\prime \nonumber \\
&=& \partial_\tau \tilde{\mathcal G}(\tau,\,s) +  \int_s^\tau \tilde{\mathcal G}(\tau,\,s^\prime)\, \partial_{s^\prime}  W(s^\prime,\,s)\, ds^\prime,
\end{eqnarray}
where we have used the fact $\partial_\tau \tilde{\mathcal G}(\tau,\,s^\prime) = -  \partial_{s^\prime} \tilde{\mathcal G}(\tau,\,s^\prime)$ following from the construction of the parametrix $\tilde{\mathcal G}(\tau,\,s^\prime)$. By the construction of $W(\tau,\,s)$, we can easily see that
\begin{equation}\label{N2}
\int_s^\tau \tilde{\mathcal G}(\tau,\,s^\prime)\, \partial_{s^\prime} W(s^\prime,\,s)\,U_0\, ds^\prime \in L^2((s,\,T);\,L^2(D)).
\end{equation}
In addition, using the expression of the parametrix $\tilde{\mathcal G}(\tau,\,s^\prime)$ and moving the $x$ derivatives of $\partial_x^\alpha\tilde{\mathcal G}(\tau,\,s^\prime)\, U_0$ with $\alpha\in{\mathbb Z}_+^2,\,|\alpha|\le 2$ as much as possible to $U_0$ by integration by parts, we have
\begin{equation}\label{N3}
\tilde{\mathcal G}(\tau,\,s)\, U_0 \in L^2((s,\,T);\,H^2(D)).
\end{equation}
Here we note that we need to take each $U_j$ of the open covering $\{U_j:=U(x_0^{(j)})\}_{j=1}^J$ of $\partial D$ small enough so that the above integration by parts yields \eqref{N3}. Then, using the equation for $\tilde{\mathcal G}(\tau,\,s^\prime)\, U_0$, we have
\begin{equation}\label{N4}
\partial_\tau \tilde{\mathcal G}(\tau,\,s^\prime)\, U_0 \in L^2((s^\prime,\,T);\,L^2(D)).
\end{equation}
Hence, we conclude from \eqref{N1}, \eqref{N2} and \eqref{N4} that
$\partial_\tau \mathcal G(\tau,\,s)U_0 \in L^2((s,\,T);\,L^2(D))$.

Next, let us show $\partial_x^\alpha \mathcal G(\tau,\,s)U_0 \in L^2((s,\,T);\,L^2(D))$ for $|\alpha|\le2$. Indeed, by defining $\tilde U = U -U_0$, we have
\begin{equation*}
\left(
\begin{array}{cc}
-\Delta & 0 \\ 0  &  -k\Delta
\end{array} \right) \tilde U = \left(
\begin{array}{lr}
-\Delta & 0 \\ 0  &  -k\Delta
\end{array} \right) U_0 - \partial_\tau \tilde U  \qquad \mathrm{ in }\; D\times (s,\,T)
\end{equation*}
together with the homogeneous transmission condition on $(\partial D)\times (s,\,T)$, which is an elliptic interior transmission problem. Then, by the construction of the parametrix for \eqref{eq:L_ITP-2} with $\tau=0$, there exists a parametrix for this elliptic interior transmission problem (for this readers can also see \cite{K-N-S}), and hence we have $\tilde U \in L^2((s,\,T);\,H^2(D))$. In summary, we have proven that $U \in F((s,\,T))$ for $U_0\in C_0^\infty(D)$.

\medskip
Now we can easily show \eqref{uq2} which gives the uniqueness of the Green function $\mathcal G$ for \eqref{Cauchy_0}.
Actually, it follows from \eqref{uq3}, and the expressions
\begin{equation*}
\int_D U(z,\,t)\cdot Z(z,\,t)\,dz = \int_D \left( \int_D \mathcal G(x,\,t;\,y,\,s)U_0(y)\,dy \right)\cdot Z_0(x)\,dx,
\end{equation*}
\begin{eqnarray*}
\int_D U(z,\,s)\cdot Z(z,\,s)\,dz &=& \int_D U_0(y)\cdot \left( \int_D \mathcal H(y,\,s;\,x,\,t)Z_0(x)\,dx  \right)\,dy\\
&=& \int_D \left( \int_D \mathcal H^*(y,\,s;\,x,\,t)U_0(y)\,dy  \right) \cdot Z_0(x)\,dx.
\end{eqnarray*}
Then, by recalling \eqref{G}, we have thus shown that the interior transmission problem \eqref{eq:ITP-2} also has a unique Green function.

Finally, we show the uniqueness of solutions to \eqref{eq:ITP-2} in $W((0,\,T))$. Again we use a duality argument. Let $x\in D$, $0<\tau<t<T$, $N\in C_0^\infty(D_T)$ and $U(\tau)=\mathbb {G}(\tau)N$.  Likewise the previous duality argument to prove the uniqueness of the Green function, we consider
\begin{equation}
\Psi(\tau)=\int_D \mathcal{H}(y,\,\tau;\,x,\,t)\, U(y,\,\tau)\,dy.
\end{equation}
Here note that due to $\tau<t$, $\mathcal{H}(y,\,\tau;\, x,\,t)$ is smooth enough with respect to $(y,\,\tau)$ (see Remark \ref{Arima_R}). Then by a similar argument leading to \eqref{uq3}, we have
\begin{equation}
\frac{d}{d\tau}\Psi(\tau)=\int_D \mathcal{H}(y,\,\tau;\,x,\,t)\,N(y,\,\tau)\,dy.
\end{equation}
Integrating this from $0$ to $t$ with respect to $\tau$, we have
\begin{equation}\label{representation}
U(x,\,t)=\int_0^t\int_D\,\mathcal{H}(y,\,\tau;\, x,\,t)\, N(y,\,\tau)\, dy\,d\tau.
\end{equation}
Suppose there are two solutions $U(t)=\mathbb {G}(\tau)N$ and $U^\prime(t)=\mathbb {G}(\tau)N$ for given $N\in L^2(D_T)$. Let $N_j\in C_0^\infty(D_T),\,j\in\mathbb{N}$ be a sequence which approximates $N$ in $L^2(D_T)$. Then from \eqref{representation}, we have
\begin{eqnarray*}
U_j^\prime(x,\,t)=U_j(x,\,t),\qquad (x,\,t)\in D_T,\;j\in\mathbb{N}.
\end{eqnarray*}
Since $\mathbb{H}: L^2(D_T)\rightarrow L^2((0,\,T);\, H^1(D))$ is bounded, we have $U^\prime(t)=U(t),\,t\in(0,\,T)$, where $\mathbb{H}:\,L^2(D_T)\to L^2((0,\,T);\, H^1(D))$ is the operator with the kernel $\mathcal{H}(y,\,s;\,x,\,t)$ which is also
bounded. This completes the proof of the uniqueness of solutions to \eqref{eq:ITP-2} in $W((0,\,T))$.

\section{Concluding remarks}\label{conclusion}
\setcounter{equation}{0}

In this paper, we investigated the interior transmission problem for the diffusion equation which is a non-classical initial boundary value problem for a pair of the diffusion equations with coupled boundary conditions. This work stemmed from our previous studies on the sampling method for reconstructing unknown inclusions in a diffusive conductor from boundary measurements. The unique solvability of the interior transmission problem was obtained by showing the existence and uniqueness of its Green function. Our approach of constructing the Green function is based on the theory of pseudo-differential operators with a large parameter. We adapted Seeley's argument of analyzing boundary value problems for elliptic equations. There are three important facts we used for the construction. First, the Lopatinskii matrix for the ITP in the Laplace domain is invertible, which enables us to construct a parametrix for the ITP in the Laplace domain. Second, the amplitude of this parametrix satisfies the assumptions of Lemma \ref{Arima} so that we can have a parametrix for the parabolic ITP by the inverse Laplace-Fourier transform. Third, the Levi method was used to compensate this parametrix to generate the Green function for the parabolic ITP, which is the advantage of considering ITP for parabolic equations.

Apart from Seeley's argument, there is an another argument to construct the Green function using the theory of pseudo-differential operators with a large parameter and boundary layer potentials; see the argument given in \cite{Greiner1971}. The principal part of our Green function is explicitly given and it can be efficiently used to analyze the asymptotic behavior of indicator function of the linear sampling method for inverse boundary problems such as active thermography and optical thermography. We assumed for simplicity that the heat conductivities in $D$ for $u$ and $v$ are both homogeneous and isotropic. The generalization to the inhomogeneous and anisotropic case is almost straightforward.

\bigskip
\noindent{\bf Acknowledgements:} The authors would like to thank the referee for his careful reading and valuable suggestions, which made the paper much improved. This work is supported by National Natural Science Foundation of China (No. 11671082) and Grant-in-Aid for Scientific Research (15K21766 and 15H05740) of the Japan Society for the Promotion of Science. The second author is also sponsored by Qing Lan Project of Jiangsu Province.

\section*{Appendix}\label{Appendix1}

\subsection*{1 Solving the linear system for $L=2$ in Theorem \ref{param}}
\setcounter{equation}{0}
\renewcommand{\theequation}{A.\arabic{equation}}

In this appendix, we show some details of the proof for $L=2$ in Theorem \ref{param}. From \eqref{L2_5}--\eqref{L2_7}, we obtain the following system of equations for constants $C_j\,(j=1,\,\cdots,\,12)$ involved in \eqref{L2_1}--\eqref{L2_4}:
\begin{equation}\label{A2_1}
(2-\ell)(C_1 - C_4)\exp\left((\lambda_{+}-\lambda_{-}) y_3\right) + (\ell -1) (C_2 - C_5) \exp\left((\lambda_{+} - \mu_{-}) y_3\right) + (2-\ell)(C_3 - C_6) =0,
\end{equation}
\begin{eqnarray}\label{A2_2}
&&\lambda_{+}(2-\ell) (C_1 - C_4) \exp\left((\lambda_{+} - \lambda_{-}) y_3\right) \nonumber\\ &&  +\lambda_{+}(\ell-1) (C_2 - C_5) \exp\left((\lambda_{+} - \mu_{-}) y_3\right)  + (2-\ell)(\lambda_{-} C_3 - \lambda_{+} C_6)  = - (2-\ell) (F_{1,3} - F_{1,4}),
\end{eqnarray}
\begin{equation}\label{A2_3}
(2-\ell)(C_7 - C_{10}) \exp\left((\mu_{+} - \lambda_{-}) y_3\right) + (\ell-1)(C_8 - C_{11})\exp\left((\mu_{+} - \mu_{-}) y_3\right) + (\ell-1)( C_9 - C_{12}) = 0,
\end{equation}
\begin{eqnarray}\label{A2_4}
&&\mu_{+} (2-\ell) (C_7 - C_{10}) \exp\left((\mu_{+} - \lambda_{-}) y_3\right) \nonumber \\
&& + \mu_{+}(\ell-1)(C_8 - C_{11}) \exp\left((\mu_{+} - \mu_{-}) y_3\right) + (\ell-1)( \mu_{-} C_9 - \mu_{+} C_{12}) = - (\ell -1)(F_{1,7} - F_{1,8}),
\end{eqnarray}
\begin{eqnarray}
\label{A2_5}
&& (2-\ell)(C_1 + C_3 - C_7) \exp\left(-\lambda_{-} y_3\right) + (\ell-1)(C_2 - C_8 - C_9)\exp\left(-\mu_{-} y_3\right) \nonumber\\
& = & - (2-\ell)\displaystyle\sum_{l=1}^2 (F_{l,1}+F_{l,3}-F_{l,5})(-y_3)^l \exp\left(-\lambda_{-} y_3\right) - (\ell-1)\displaystyle \sum_{l=1}^2
(F_{l,2}-F_{l,6}-F_{l,7})(-y_3)^l \exp\left(-\mu_{-} y_3\right) \nonumber\\
& =:&(2-\ell) A_3 \exp\left(-\lambda_{-} y_3\right) + (\ell-1) B_3 \exp\left(-\mu_{-} y_3\right),
\end{eqnarray}
\begin{eqnarray*}
&&(i\displaystyle\sum_{j=1}^2m_{3j}^0\xi_j + \lambda_{+} m_{33}^0)((2-\ell)C_1 \exp\left(-\lambda_{-} y_3\right) + (\ell-1)C_2 \exp\left(-\mu_{-} y_3\right)) \nonumber\\
&&+(i\sum_{j=1}^2m_{3j}^0\xi_j + \lambda_{-}m_{33}^0)(2-\ell)C_3 \exp\left(-\lambda_{-} y_3\right)\nonumber\\
&&- k(i\displaystyle\sum_{j=1}^2m_{3j}^0\xi_j + \mu_{+}m_{33}^0)((2-\ell)C_7 \exp\left(-\lambda_{-} y_3\right) + (\ell-1) C_8 \exp\left(-\mu_{-} y_3\right)) \nonumber\\
&& - k(i\sum_{j=1}^2m_{3j}^0\xi_j + \mu_{-} m_{33}^0)(\ell-1) C_9 \exp\left(-\mu_{-}y_3\right)\nonumber
\\
& = & -(i\displaystyle\sum_{j=1}^2m_{3j}^0\xi_j)[\sum_{l=1}^2(2-\ell)(F_{l,1}+F_{l,3})(-y_3)^l
\exp\left(-\lambda_{-} y_3\right) + \sum_{l=1}^2(\ell-1) F_{l,2}(-y_3)^l \exp\left(-\mu_{-}y_3\right)]\nonumber\\
&& -m_{33}^0\big[(2-\ell)\displaystyle\{\sum_{l=1}^2l(F_{l,1}+F_{l,3})(-y_3)^{l-1}+\sum_{l=1}^2(\lambda_{+}F_{l,1}
+ \lambda_{-} F_{l,3})(-y_3)^l\} \exp\left(-\lambda_{-} y_3\right)\nonumber\\
\end{eqnarray*}
\begin{eqnarray}\label{A2_6}
&& + (\ell-1)\displaystyle\{\sum_{l=1}^2lF_{l,2}(-y_3)^{l-1}+\lambda_{+} \sum_{l=1}^2 F_{l,2}(-y_3)^l\} \exp\left(-\mu_{-} y_3\right)\big]\nonumber\\
&& +(ik\displaystyle\sum_{j=1}^2m_{3j}^0\xi_j)[\sum_{l=1}^2(2-\ell) F_{l,5}(-y_3)^l \exp\left(-\lambda_{-} y_3\right) + \sum_{l=1}^2 (\ell-1)(F_{l,6}+F_{l,7})(-y_3)^l \exp\left(-\mu_{-} y_3\right)]\nonumber\\
&& +km_{33}^0\big[(2-\ell)\{\displaystyle\sum_{l=1}^2lF_{l,5}(-y_3)^{l-1}+\mu_{+}\sum_{l=1}^2F_{l,5} (-y_3)^l\} \exp\left(-\lambda_{-} y_3\right)\nonumber\\
&& +(\ell-1)\{\displaystyle\sum_{l=1}^2l(F_{l,6}+F_{l,7})(-y_3)^{l-1}+\sum_{l=1}^2(\mu_{+}F_{l,6}+ \mu_{-} F_{l,7})(-y_3)^l\} \exp\left(-\mu_{-} y_3\right)\big]\nonumber\\
& =:& (2-\ell)A_4 \exp\left(-\lambda_{-} y_3\right) + (\ell-1)B_4 \exp\left(-\mu_{-}y_3\right),
\end{eqnarray}
where $\mathrm{ord}\,A_3=\mathrm{ord}\,B_3=-2$ and $\mathrm{ord}\, A_4=\mathrm{ord}\, B_4=-1$.

\medskip
We remark that, looking at the structures of the amplitudes, $C_3,\,C_9$ and the four linear combinations given by $(2-\ell)C_1 \exp\left(-\lambda_- y_3\right) + (\ell-1)C_2 \exp\left(-\mu_- y_3\right)$, $(2-\ell)C_4 \exp\left(-\lambda_- y_3\right)+ (\ell-1)C_5 \exp\left(-\mu_- y_3\right) + (2-\ell)C_6 \exp\left(-\lambda_+ y_3\right)$, $(2-\ell)C_7\exp\left(-\lambda_ - y_3\right) + (\ell-1)C_8 \exp\left(-\mu_- y_3\right)$, $(2-\ell)C_{10} \exp\left(-\lambda_- y_3\right) + (\ell-1)C_{11}\exp\left(-\mu_- y_3\right) + (\ell-1)C_{12} \exp\left(-\mu_+ y_3\right)$ are the substantial unknowns for the system \eqref{A2_1}--\eqref{A2_6}. Hence, $\mathbf a_{\ell,-2},\, \mathbf b_{\ell,-2},\, \mathbf d_{\ell,-2},\, \mathbf e_{\ell,-2}$ are uniquely determined from \eqref{A2_1}--\eqref{A2_6}. We show the details as follows.

From \eqref{A2_1}--\eqref{A2_4}, we can easily get
\begin{equation}\label{second_0}
C_3=\frac{F_{1,3}-F_{1,4}}{\lambda_{+}-\lambda_{-}},\quad C_9=\frac{F_{1,7}-F_{1,8}}{\mu_{+}-\mu_{-}},
\end{equation}
where $\mathrm{ord}\, C_3=\mathrm{ord}\, C_9=-2$. Then \eqref{A2_5} and \eqref{A2_6} become
\begin{eqnarray}\label{second_1}
&&(2-\ell)C_1 \exp\left(-\lambda_{-} y_3\right) + (\ell-1)C_2 \exp\left(-\mu_{-} y_3\right) - (2-\ell)C_7
\exp\left(-\lambda_{-} y_3\right) - (\ell-1)C_8 \exp\left(-\mu_{-} y_3\right) \nonumber\\
& = & (2-\ell) (A_3 - C_3) \exp\left(-\lambda_{-} y_3\right) + (\ell-1)(B_3 + C_9) \exp\left(-\mu_{-} y_3\right) \nonumber\\
& =: & (2-\ell) A_5 \exp\left(-\lambda_{-} y_3\right) + (\ell-1) B_5 \exp\left(-\mu_{-} y_3\right)
\end{eqnarray}
and
\begin{eqnarray}\label{second_2}
&&(i\displaystyle\sum_{j=1}^2m_{3j}^0\xi_j + \lambda_{+} m_{33}^0)( (2-\ell)C_1 \exp\left(-\lambda_{-} y_3\right) + (\ell-1)C_2 \exp\left(-\mu_{-} y_3\right))\nonumber \\
&&-k(i\displaystyle\sum_{j=1}^2m_{3j}^0\xi_j + \mu_{+} m_{33}^0)((2-\ell)C_7 \exp\left(-\lambda_{-} y_3\right) + (\ell-1)C_8 \exp\left(-\mu_{-} y_3\right)) \nonumber\\
&=&(2-\ell)\big\{A_4-(i\displaystyle\sum_{j=1}^2 m_{3j}^0\xi_j + \lambda_{-} m_{33}^0)C_3\big\} \exp\left(-\lambda_{-} y_3\right) \nonumber\\
&&+(\ell-1)\big\{B_4+k(i\displaystyle\sum_{j=1}^2 m_{3j}^0\xi_j + \mu_{-} m_{33}^0)C_9\big\}\exp\left(-\mu_{-} y_3\right) \nonumber\\
&=:& (2-\ell) A_6 \exp\left(-\lambda_{-} y_3\right) + (\ell-1) B_6 \exp\left(-\mu_{-} y_3\right),
\end{eqnarray}
respectively, where $\mathrm{ord}\, A_5=\mathrm{ord}\, B_5=-2$ and $\mathrm{ord}\, A_6=\mathrm{ord}\,B_6=-1$. By direct calculations, we derive from \eqref{second_1} and \eqref{second_2} that
\begin{eqnarray}\label{second_3}
&&(2-\ell)C_7 \exp\left(-\lambda_{-} y_3\right) + (\ell-1)C_8 \exp\left(-\mu_{-} y_3\right) \nonumber\\
& = & \Big\{k(i\sum_{j=1}^2m_{3j}^0\xi_j + \mu_{+}
m_{33}^0)-(i\sum_{j=1}^2m_{3j}^0\xi_j + \lambda_{+}
m_{33}^0)\Big\}^{-1}\nonumber\\
& & \times \Big[(2-\ell)\Big\{(i\sum_{j=1}^2 m_{3j}^0\xi_j +
\lambda_{+} m_{33}^0)A_5 - A_6\Big\}\exp\left(-\lambda_{-} y_3\right) \nonumber\\
& & \qquad +(\ell-1)\Big\{(i\sum_{j=1}^2m_{3j}^0 \xi_j +
\lambda_{+}m_{33}^0)B_5-B_6\Big\}\exp\left(-\mu_{-} y_3\right)\Big] \nonumber\\
& =:& (2-\ell) A_7 \exp\left(-\lambda_{-} y_3\right) + (\ell-1) B_7 \exp\left(-\mu_{-} y_3\right),
\end{eqnarray}
\begin{eqnarray}\label{second_4}
&&(2-\ell)C_1 \exp\left(-\lambda_{-} y_3\right) + (\ell-1) C_2 \exp\left(-\mu_{-} y_3\right)\nonumber\\ & = & (2-\ell)(C_7 + A_5)
\exp\left(-\lambda_{-} y_3\right) + (\ell-1)(C_8 + B_5) \exp\left(-\mu_{-} y_3\right)\nonumber\\
& = & (2-\ell)(A_5 + A_7) \exp\left(-\lambda_{-} y_3\right) + (\ell-1)(B_5 + B_7) \exp\left(-\mu_{-}y_3\right)
\end{eqnarray}
with $\mathrm{ord}\,A_7=\mathrm{ord}\,B_7=-2$. Consequently, we have
\begin{eqnarray}\label{second_5}
&& (2-\ell)C_4 \exp\left(-\lambda_{-} y_3\right) + (\ell-1)C_5 \exp\left(-\mu_{-} y_3\right) + (2-\ell)C_6
\exp\left(-\lambda_{+} y_3\right) \nonumber\\
& = & \exp\left(-\lambda_{+} y_3\right) [(2-\ell) C_4 \exp\left((\lambda_{+} - \lambda_{-}) y_3\right) +
(\ell-1)C_5 \exp\left((\lambda_{+} - \mu_{-}) y_3\right) + (2-\ell) C_6]\nonumber \\
& = & \exp\left(-\lambda_{+} y_3\right) [(2-\ell) C_1 \exp\left((\lambda_{+} - \lambda_{-}) y_3\right) +
(\ell-1)C_2 \exp\left((\lambda_{+} - \mu_{-}) y_3\right) + (2-\ell)C_3] \nonumber\\
& = & (2-\ell)(A_5 + A_7) \exp\left( -\lambda_{-} y_3\right) + (\ell-1) (B_5 + B_7) \exp\left(-\mu_{-} y_3\right) \nonumber\\
& & + (2-\ell)\frac{F_{1,3}-F_{1,4}}{\lambda_{+} - \lambda_{-}} \exp\left(-\lambda_{+} y_3\right),
\end{eqnarray}
\begin{eqnarray}\label{second_6}
&&(2-\ell)C_{10} \exp\left(-\lambda_{-} y_3\right) + (\ell-1)C_{11} \exp\left( -\mu_{-}
y_3\right) + (\ell-1)C_{12} \exp\left(-\mu_{+} y_3\right) \nonumber\\
& = & \exp\left(-\mu_{+} y_3\right) [(2-\ell)C_{10} \exp\left((\mu_{+} - \lambda_{-}) y_3\right) +
(\ell-1)C_{11} \exp\left((\mu_{+} - \mu_{-}) y_3\right) + (\ell-1)C_{12}] \nonumber\\
& = & \exp\left(-\mu_{+} y_3\right)[(2-\ell) C_7 \exp\left((\mu_{+} - \lambda_{-}) y_3\right) +
(\ell-1)C_8 \exp\left((\mu_{+} - \mu_{-}) y_3\right) + (\ell-1)C_9] \nonumber\\
& = & (2-\ell)A_7 \exp\left(-\lambda_{-} y_3\right) + (\ell-1)B_7 \exp\left(-\mu_{-} y_3\right) +
(\ell-1)\frac{F_{1,7}-F_{1,8}}{\mu_{+} - \mu_{-}}\exp\left(-\mu_{+} y_3\right).
\end{eqnarray}
By substituting \eqref{second_0}, \eqref{second_3}--\eqref{second_6} into \eqref{L2_1}--\eqref{L2_4}, we finally obtain the expressions of the amplitudes $\mathbf a_{\ell,-2},\, \mathbf b_{\ell,-2},\, \mathbf d_{\ell,-2}$ and $\mathbf e_{\ell,-2}$ as follows:
\begin{eqnarray*}
\mathbf a_{\ell,-2}&=& (2-\ell) \Big\{\sum_{l=1}^{2} F_{l,1}(x_3-y_3)^l +
(A_5+A_7)\Big\}\exp\left(\lambda_{+} x_3 - \lambda_{-} y_3 - \tau s -iy^\prime\cdot\xi^\prime\right)
\\
& & + (\ell-1) \Big\{\sum_{l=1}^{2}
F_{l,2}(x_3-y_3)^l + (B_5+B_7)\Big\}\exp\left(\lambda_{+} x_3 - \mu_{-} y_3- \tau s -iy^\prime\cdot\xi^\prime\right)\\
& & + (2-\ell)\Big\{\sum_{l=1}^{2}
F_{l,3}(x_3-y_3)^l + \frac{F_{1,3}-F_{1,4}}{\lambda_{+} - \lambda_{-}}\Big\} \exp\left(\lambda_{-} x_3 - \lambda_{-} y_3- \tau s -iy^\prime\cdot\xi^\prime\right),
\end{eqnarray*}
\begin{eqnarray*}
\mathbf b_{\ell,-2}&=&(2-\ell) \Big\{\sum_{l=1}^{2} F_{l,1}(x_3-y_3)^l + (A_5+
A_7)\Big\}\exp\left(\lambda_{+} x_3 - \lambda_{-} y_3- \tau s -iy^\prime\cdot\xi^\prime\right)
\\
& & +(\ell-1)\Big\{\sum_{l=1}^{2}
F_{l,2}(x_3-y_3)^l + (B_5 + B_7)\Big\} \exp\left(\lambda_{+} x_3 - \mu_{-} y_3- \tau s -iy^\prime\cdot\xi^\prime\right)\\
& & +(2-\ell)\Big\{\sum_{l=1}^{2}
F_{l,4}(x_3-y_3)^l + \frac{F_{1,3}-F_{1,4}}{\lambda_{+} - \lambda_{-}}\Big\} \exp\left(\lambda_{+} x_3 - \lambda_{+} y_3- \tau s -iy^\prime\cdot\xi^\prime\right),
\end{eqnarray*}
\begin{eqnarray*}
\mathbf d_{\ell,-2}&=&(2-\ell)\Big\{\sum_{l=1}^{2} F_{l,5}(x_3-y_3)^l +
A_7\Big\}\exp\left(\mu_{+} x_3 - \lambda_{-} y_3- \tau s -iy^\prime\cdot\xi^\prime\right) \\
& & + (\ell-1)\Big\{\sum_{l=1}^{2}
F_{l,6}(x_3-y_3)^l + B_7\Big\}\exp\left(\mu_{+} x_3 - \mu_{-} y_3- \tau s -iy^\prime\cdot\xi^\prime\right)\\
& & + (\ell-1)\Big\{\sum_{l=1}^{2}
F_{l,7}(x_3-y_3)^l + \frac{F_{1,7}-F_{1,8}}{\mu_{+} - \mu_{-}}\Big\}\exp\left(\mu_{-} x_3 - \mu_{-} y_3- \tau s -iy^\prime\cdot\xi^\prime\right),
\end{eqnarray*}
\begin{eqnarray*}
\mathbf e_{\ell,-2}&=& (2-\ell)\Big\{\sum_{l=1}^{2} F_{l,5}(x_3-y_3)^l +
A_7\Big\}\exp\left(\mu_{+} x_3 - \lambda_{-} y_3- \tau s -iy^\prime\cdot\xi^\prime\right) \\
& & + (\ell-1)\Big\{\sum_{l=1}^{2}
F_{l,6}(x_3-y_3)^l + B_7\Big\}\exp\left(\mu_{+} x_3 - \mu_{-} y_3- \tau s -iy^\prime\cdot\xi^\prime\right)\\
& & + (\ell-1)\Big\{\sum_{l=1}^{2} F_{l,8}(x_3-y_3)^l +
\frac{F_{1,7}-F_{1,8}}{\mu_{+} - \mu_{-}}\Big\}\exp\left(\mu_{+} x_3 -
\mu_{+} y_3- \tau s -iy^\prime\cdot\xi^\prime\right).
\end{eqnarray*}
These are the forms of \eqref{amplitude1}--\eqref{amplitude2} for $L=2$.

\subsection*{2 An alternative direct proof of $\tilde{\mathcal G}(s,\,s)=I$}
\setcounter{equation}{0}
\renewcommand{\theequation}{B.\arabic{equation}}

Note that we only need to show the property $\tilde{\mathcal G}(s,\,s)=I$ for the principal part of
$\tilde{\mathcal G}(t,\,s)$, because
the remaining part vanishes at $t=s$ which can be easily seen by using Lemma \ref{Arima}. We will only show this property for $\tilde{\mathcal{G}}^\prime_1(t,\,s)$ defined by the restriction of $\tilde{G}^\prime_1(t,\,s)=\tilde{G}^\prime_1(s)$ to $t\ge s$. Recall the definition of parametrix $\tilde{G}^\prime(t,\,s)$ with the first column vector $(\tilde{G}_1^\prime(t,\,s),\,\tilde{H}_1^\prime(t,\,s))$ defined by \eqref{parametrices tilde G_1 prime} and the second column vector
$(\tilde{G}_2^\prime(t,\,s),\,\tilde{H}_2^\prime(t,\,s))$ defined in a similar way as $(\tilde{G}_1^\prime(t,\,s),\,\tilde{H}_1^\prime(t,\,s))$. The $j=0$ term of \eqref{parametrices tilde G_1 prime} restricted to $t\ge s$ obviously satisfies this property on $\{(x,\,y)\in\textrm{supp} \,\varphi\times\textrm{supp}\,\varphi\}$. Then by Remark \ref{estimate of principal part of parametrix}, we only need to show the property for $\tilde{\dot{\mathcal{G}}}^+_{1(2)}(t,\,s)$ defined by the restriction of ${\dot G}^+_{1(2)}(x,\,t;\,y,\,s)$ (see \eqref{dot G^+_1(29)}) to $t\ge s$. By Remark \ref{symbols in Laplace domain} (ii), for any $f(y)\in C_0^\infty({\mathbb R}_-^3)$, $t>s$,
\begin{eqnarray}
(\tilde{\dot{\mathcal{G}}}^+_{1(2)}(t,\,s)f)(x,\,t)&=&\frac{-i}{(2\pi)^{4}}\displaystyle\int_{{\mathbb R}_y^3\times{\mathbb R}_\xi^3}\int_{\sigma-i\infty}^{\sigma+i\infty}\exp\big((t-s)\tau+i(x-y)\cdot\xi\big)\mathcal{J}^{-1}(y)
(\tau+(M^1\xi)\cdot\xi)^{-1}\,f(y)\,d\tau\,dy\,d\xi \nonumber\\
&=&\frac{1}{(2\pi)^{3}}\displaystyle\int_{{\mathbb R}_y^3\times{\mathbb R}_\xi^3}\exp\big(i(x-y)\cdot\xi-(t-s)((M^1\xi)\cdot\xi)\big)\mathcal{J}^{-1}(y)\,f(y)\,d\xi\,dy.
\end{eqnarray}
Hence, $\lim_{t\downarrow s}\,\tilde{\dot{\mathcal{G}}}^+_{1(2)}(t,\,s)(x,\,y)=\mathcal{J}^{-1}(y)\delta(x-y)I$ in terms of coordinates introduced by the coordinates transformation $\Phi_{x_0}$ flattening $\partial D$ near $x_0\in\partial D$. This completes proving the property.

\subsection*{3 A fact for the $L^2$ boundedness}
\setcounter{equation}{0}
\renewcommand{\theequation}{C.\arabic{equation}}
Here we give a well known fact which we used to show that the Green function $\mathbb{G}$ for \eqref{eq:ITP-2} and the Green function $\mathbb{H}$ for the adjoint problem of \eqref{eq:ITP-2} are bounded operators from $L^2(D_T)$ to $L^2((0,\,T);\, H^1(D))$.

Let $(X_j, \,B_j, \,\mu_j),\,\,j=1,\,2$ be sigma finite complete measured spaces and $(Y,\,B,\,\mu)$ be the completion of the product measure space of them. Suppose $K(x_1,\,x_2)$ is a $B$ measurable function such that
\begin{equation}
\left\{
\begin{array}{l}
\displaystyle\int_{X_1}|K(x_1,\,x_2)|\,d\mu_1(x_1)\le M_1,\qquad \mu_2\textrm{-a.e.}\,\,x_2\in X_2,\\
\displaystyle\int_{X_2}|K(x_1,\,x_2)|\,d\mu_2(x_2)\le M_2,\qquad \mu_1\textrm{-a.e.}\,\,x_1\in X_1.
\end{array}
\right.
\end{equation}
Then for any $f\in L^2(X_2,\,d\mu_2)$, we have
\begin{equation}
\Vert Kf\Vert_{L^2(X_1,\,d\mu_1)}\le \sqrt{M_1\,M_2}\,\Vert f\Vert_{L^2(X_2,\,d\mu_2)},
\end{equation}
where $Kf$ is given by
$$
(Kf)(x_1)=\int_{X_2} K(x_1,\,x_2)\,f(x_2)\, d\mu_2(x_2).
$$

\end{document}